\documentclass[a4paper]{article}

\usepackage{amsfonts, amsmath, amsthm, empheq, amssymb, bm, setspace, mathrsfs}
\usepackage[numbers,sort]{natbib}
\usepackage[colorlinks,pdfborder=001,linkcolor=blue,citecolor=blue]{hyperref}

\usepackage[pagewise]{lineno}

\makeatletter

\newcommand{\Rmnum}[1]{\expandafter\@slowromancap\romannumeral #1@}
\@addtoreset{equation}{section}
\makeatother

\allowdisplaybreaks[4]

\newtheorem{thm}{Theorem}[section]

\newtheorem{lem}[thm]{Lemma}

\title{\Large\bf\boldmath
Carleman estimates for a magnetohydrodynamics system and application to inverse source problems}

\author{\large Xinchi HUANG$^*$$^{1,2}$ and Masahiro YAMAMOTO$^{1,3,4}$}
\date{}

\begin{document}
\maketitle

\renewcommand{\thefootnote}{\fnsymbol{footnote}}
\footnotetext{


* Corresponding author.

$^1$ Graduate School of Mathematical Sciences, The University of Tokyo, 3-8-1 Komaba, 
Meguro, Tokyo, 153-8914, Japan.

$^2$ JSPS Postdoctoral Fellowships for research in Japan.

$^3$ Honorary Member of Academy of Romanian Scientists, 
Ilfov, nr. 3, Bucuresti, Romania.

$^4$ Correspondence member of Accademia Peloritana dei Pericolanti,
Palazzo Universit\`a, Piazza S. Pugliatti 1 98122 Messina, Italy. 


}

\begin{abstract}
\noindent 
In this article, we consider a linearized magnetohydrodynamics system
for incompressible flow in a three-dimensional bounded domain. 
We first prove two kinds of Carleman estimates.  
This is done by combining the Carleman estimates for the parabolic and the elliptic equations.  
Then we apply the Carleman estimates to prove H\"older type stability results 
for some inverse source problems. 
\end{abstract}

\section{Introduction}
\label{sec-intro}

The magnetohydrodynamics equations (the MHD equations in short)
are concerned with the magnetic properties of electrically conducting fluids 
such as plasmas, liquid metals and salt water.

Before the formulation, we introduce the following notations. 
Let $\Omega\subset\mathbb{R}^3$ be a bounded domain with smooth boundary, 
and let $n=n(x)$ be the outward unit normal vector to $\partial\Omega$ at
$x$, and $\partial_n\varphi = \nabla \varphi \cdot n$. 
Moreover $\cdot^T$ denotes the transpose of matrices or vectors, and let 
$\partial_t = \frac{\partial}{\partial t}, \; 
\partial_j = \frac{\partial}{\partial x_j}, j = 1,2,3, \; 
\Delta = \sum_{j=1}^3 \partial_j^2, \; 
\nabla = (\partial_1, \partial_2, \partial_3)^T, \; 
\nabla_{x,t} = (\nabla,\partial_t)^T $. 
We also use the following notations:
$$
(w \cdot \nabla) v := 
\left(\sum_{j=1}^3 w_j\partial_j v_1, 
\sum_{j=1}^3 w_j\partial_j v_2, 
\sum_{j=1}^3 w_j\partial_j v_3 \right)^T,
$$
$$
v\times w := 
(v_2w_3 - v_3w_2, 
v_3w_1 - v_1w_3, 
v_1w_2 - v_2w_1)^T,
$$
$$
\mathrm{div}\, w := \sum_{j=1}^3 \partial_j w_j, \qquad
\mathrm{rot}\, w := 
(\partial_2w_3 - \partial_3w_2, 
\partial_3w_1 - \partial_1w_3, 
\partial_1w_2 - \partial_2w_1)^T
$$
for vectors $v = (v_1,v_2,v_3)^T$, $w = (w_1,w_2,w_3)^T$.

By using the divergence-free condition of the magnetic field, 
we simplify Maxwell's equations into a second-order parabolic equation. 
Then for the velocity field $u=(u_1,u_2,u_3)^T$, the pressure $p$ and 
the magnetic field intensity $H = (H_1,H_2,H_3)^T$, 
the MHD equations are described as follows:  
$$
\left\{
\begin{aligned}
& \ \partial_t u - \nu\Delta u + (u\cdot\nabla) u 
- \mu\, \mathrm{rot}\, H\times H + \nabla p = F, \\
& \ \partial_t H - \frac{1}{\sigma\mu}\Delta H - \mathrm{rot}(u\times H) = 0, \\
& \ \mathrm{div}\, u = 0, \quad \mathrm{div}\, H = 0.
\end{aligned}   
\right.                                                 
$$
Here $\nu$ describes the viscosity of the fluids, while $\sigma$ and $\mu$ 
denote the electric conductivity and the magnetic permeability, respectively. 
We assume that they are all positive constants.
For more details about the derivation of the MHD equations, 
we refer to Li and Qin \cite{LQ13} for example. 
Moreover, $\mathrm{div}\, H = 0$ is used for the proof of one of 
our Carleman estimates (Theorem \ref{thm:CE2}).

There are several papers concerning the direct problems of the MHD equations. 
Since our main purpose of this article is the inverse problems, 
we do not explain more on the direct problems and 
we refer to Ladyzhenskaya and Solonnikov \cite{LS78}, as a classical work, 
Li \cite{L18} and the references therein. 

The purpose of this article is the stability for inverse source problems, 
more precisely, the determination of the spatially varying factor in the source term. 
Our main strategy is based on a method by Bukhgeim-Klibanov \cite{BK81} and 
Klibanov \cite{K92}, applying so-called Carleman estimates. 

Carleman estimate is an $L^2$-weighted estimate for the solution to 
a partial differential equation with large parameter, 
which is nowadays known as a powerful tool in treating 
the uniqueness and stability estimates for inverse problems. 
Actually, the Carleman estimate itself and its applications to inverse problems 
have been intensively studied for a variety of partial differential equations 
in mathematical physics including, for example,
\begin{itemize}
\item Transport equations: 
Cannarsa, Floridia and Yamamoto \cite{CFY19}, 
Cannarsa, Floridia, G\"olgeleyen and Yamamoto \cite{CFGY19};

\item Hyperbolic equations: 
Bellassoued and Yamamoto \cite{BY17}, 
Imanuvilov and Yamamoto \cite{IY01};

\item Parabolic equations: 
Imanuvilov and Yamamoto \cite{IY98}, 
Yamamoto \cite{Y09};

\item Lam\'{e} equation and the Navier-Stokes equations: 
Bellassoued, Imanuvilov and Yamamoto \cite{BIY08, BIY16}, 
Choulli, Imanuvilov, Puel and Yamamoto \cite{CIPY13}.
\end{itemize} 

To the authors' best knowledge, there are few papers on Carleman estimates 
for the MHD equations. 
For the MHD equations, Hav\^arneanu, Popa and Sritharan \cite{HPS06,HPS07} 
proved a Carleman estimate and established the exact controllability. 
In their Carleman estimate, the observation of the first-order spatial derivative 
of the source $F$ is necessary, and so it is not suitably designed for inverse problems. 
Besides, Huang \cite{H21} established Carleman estimates for the MHD equations 
and investigate inverse coefficient problems by the measurements of velocity field, 
magnetic field and pressure on some sub-boundary. 
However, some additional data on the derivative of pressure are imposed due to 
the technical proof. 
In this article, we apply the recent idea of Huang, Imanuvilov, Yamamoto \cite{HIY20} 
and Imanuvilov, Lorenzi, Yamamoto \cite{ILY21} to eliminate the additional data and 
propose simplified proofs for the Carleman estimates as well as the stability estimates 
for the inverse source problems. 
Furthermore, under some assumptions on the source term 
we discuss also the inverse source problem without the knowledge of the pressure, 
that is, inverse source problem 2 in the following context. 

Now we formulate our inverse problems and state the main results. 
For simplicity, we set $\kappa = \frac{1}{\sigma\mu} > 0$ 
and we consider the following linearized MHD equations:
\begin{equation}
\label{sy:MHD}
\left\{
\begin{aligned}
&\ \partial_t u - \nu\Delta u + (A^{(1)}\cdot\nabla) u 
+ (u\cdot\nabla) A^{(2)} + L_1H + \nabla p = F(x,t)   &&\quad \mbox{in} \ \Omega\times (0,T), \\
&\ \partial_t H - \kappa\Delta H + (A^{(3)}\cdot\nabla) H 
+ (H\cdot\nabla) A^{(4)} + L_2u = 0                          &&\quad \mbox{in} \ \Omega\times (0,T), \\
&\ \mathrm{div}\, u = 0, \quad \mathrm{div}\, H = 0  &&\quad \mbox{in} \ \Omega\times (0,T),
\end{aligned}
\right.                
\end{equation}
where the coupling operators $L_1, L_2$ are defined by 
\begin{equation}
\label{def:L1L2}
\left\{
\begin{aligned}
&\ L_1 H := (C^{(1)}\cdot\nabla) H + (H\cdot\nabla) C^{(2)} + \nabla (C^{(3)}\cdot H), \\
&\ L_2 u := (u\cdot\nabla) C^{(4)} + (C^{(5)}\cdot \nabla)u, 
\end{aligned}
\right.
\end{equation} 
and the vector coefficients $A^{(j)}, C^{(k)}$, $j=1,2,3,4$, $k=1,2,3,4,5$ 
are supposed to be sufficiently smooth. 
Henceforth we arbitrarily fix $0<t_0<T$, $0<\delta<\min\{t_0,T-t_0\}$ 
and we denote
$$
I := (t_0-\delta,t_0+\delta), \quad Q := \Omega \times I. 
$$
Moreover, let $H^k(\Omega)$, etc., denote usual Sobolev spaces 
(e.g., Adams and Fournier \cite{AF03}). 
Let $\gamma_0\in (\Bbb N \cup \{0\})$, 
$\gamma = (\gamma_1,\gamma_2,\gamma_3) \in (\Bbb N \cup \{0\})^3$, 
$\partial_x^{\gamma} = \partial_1^{\gamma_1} \partial_2^{\gamma_2} \partial_3^{\gamma_3}$ 
and $\vert \gamma \vert = \gamma_1 + \gamma_2 +\gamma_3$. 
Then we set 
\begin{equation*}
\left\{
\begin{aligned}
&\ W^{k,\infty}(D) :=
\{w; \; \partial_t^{\gamma_0} w,\partial_x^{\gamma} w \in L^\infty (D), \; 
\gamma_0\le k, |\gamma|\le k \}, 
\ k\in \Bbb{N} \\
&\ H^{k,\ell}(D) :=
\{w; \; \partial_t^{\gamma_0} w, \partial_x^{\gamma} w \in L^2(D), \; 
\gamma_0\le \ell, |\gamma|\le k \}, 
\ k,\ell\in \mathbb{N}\cup \{ 0 \}
\end{aligned}
\right.
\end{equation*}
for any sub-domain $D \subset Q$.
If there is no confusion, we do not distinguish $L^2(\Omega)$ 
with $L^2(\Omega; \mathbb{R}^3)$ and we use $L^2(\Omega)$ for both meanings. 

Henceforth, depending on different cases, we consider some of the following
conditions:
\begin{align}
\label{con:F1}
& |\partial_t^k F(x,t)| \le C|F(x,t_0)|,                                  
             \quad (x,t)\in Q,\ k=0,1. \\
\label{con:F2-1}
& |\partial_t^k \mathrm{rot}\, F(x,t)| \le C|\mathrm{rot}\, F(x,t_0)|,               
             \quad (x,t)\in Q,\ k=0,1. \\
\label{con:F2-2a}
& |\partial_t^k \mathrm{rot}\, F(x,t)| \le C(|\nabla F(x,t_0)| + |F(x,t_0)|),         
             \quad (x,t)\in Q,\ k=0,1,2.\\
\label{con:F2-2b}
& |\mathrm{div}\, F(x,t_0)| = 0, 
             \quad x\in \Omega.
\end{align}

We are ready to state the different conditional stability estimates under the conditions 
\eqref{con:F1},  \eqref{con:F2-1} and \eqref{con:F2-2a}--\eqref{con:F2-2b}, respectively.  
Our first inverse source problem can be described as follows: 

\noindent $\mathbf{Inverse\ source\ problem\ 1}:$  

Let $F$ satisfy \eqref{con:F1}. 
For arbitrarily given sub-boundary $\Gamma \subset \partial\Omega$, 
and arbitrarily given sub-domain $\Omega_0 \subset \Omega \cup \Gamma$,
determine the source term $F(\cdot,t_0)$ in $\Omega_0$ by the observation data of 
$\{ (u,\nabla u, p, H, \nabla H)|_{\Gamma \times I}, \; 
u(\cdot,t_0)|_{\Omega}, \, 
p(\cdot, t_0)|_{\Omega} \}$.

More precisely, by introducing an admissible set of the unknown function $F$: 
\begin{align*}
\mathcal{F}_{M,1} = 
\Big\{ 
& F\in H^1(0,T;L^2(\Omega)); 
\|F\|_{H^1(0,T;L^2(\Omega))} \le M, 
\mbox{ \eqref{con:F1} is satisfied} 
\Big\}
\end{align*}
with arbitrarily given constant $M>0$ and the norm of the data: 
\begin{align*}
D_1 := \ 
&\|u\|_{H^2(I;L^2(\Gamma))} 
+ \|\nabla u\|_{H^1(I;L^2(\Gamma))} 
+ \|H\|_{H^2(I;L^2(\Gamma))} 
+ \|\nabla H\|_{H^1(I;L^2(\Gamma))} \\
&+ \|p\|_{H^1(I;H^{\frac12}(\Gamma))} 
+ \|u(\cdot,t_0)\|_{H^2(\Omega)} 
+ \|\nabla p(\cdot,t_0)\|_{L^2(\Omega)},
\end{align*}
we have the following stability estimate of H\"older type for the inverse source problem.
\begin{thm}
\label{thm:stab1}
Assume that $F\in \mathcal{F}_{M,1}$ and a solution to \eqref{sy:MHD} satisfies 
$(u,p,H)\in H^{2,2}(Q)\times H^{1,1}(Q)\times H^{2,2}(Q)$ with 
\begin{align*}
\|u\|_{H^{1,2}(Q)} + \|\nabla u\|_{H^{0,1}(Q)} 
+ \|H\|_{H^{1,1}(Q)} + \|\nabla H\|_{H^{0,1}(Q)} 
+ \|p\|_{H^{0,1}(Q)} \le M.
\end{align*} 
Then there exist constants $C > 0$ and $\theta\in (0,1)$ such that
\begin{align*}
\|F(\cdot,t_0)\| _{L^2(\Omega_0)} \le C (D_1^\theta + D_1).
\end{align*}
\end{thm}

\noindent {\bf Remark 1.} 
In the above theorem and the following theorems in this section, 
the constants $C$ and $\theta$ depend on $M$, $\delta$, $\Omega$, 
the coefficients and also the choices of $t_0$, $\Omega_0$ and $\Gamma$. 
In particular, as one can see from the proof, the constant $C$ goes to infinity 
and the constant $\theta$ goes to zero as $t_0$ goes to zero, 
or $\Omega_0$ tends to $\Omega$, or $\Gamma$ tends to an empty set. 
However, we mention that the arbitrary choice of $\Omega_0$ enables us 
to prove the corresponding uniqueness results of determining 
$F(\cdot,t_0)$ in $\Omega$ immediately by contradiction. 

\noindent {\bf Remark 2.} 
According to the condition \eqref{con:F1}, one can readily see that Theorem \ref{thm:stab1} 
also implies the estimate of the source term in $\Omega_0\times I$ as follows:
\begin{align*}
\|F\| _{L^2(\Omega_0\times I)} \le C (D_1^\theta + D_1).
\end{align*}

In Theorem \ref{thm:stab1}, the data on the pressure $p$ is needed. 
However, sometimes it is not easy to measure the pressure. 
Thus, we discuss also the second inverse source problem as follows:

\noindent $\mathbf{Inverse\ source\ problem\ 2}:$  

Let $F$ satisfy \eqref{con:F2-2a} and \eqref{con:F2-2b}. 
For arbitrarily given sub-boundary $\Gamma \subset \partial\Omega$, 
and arbitrarily given sub-domain $\Omega_0 \subset \Omega \cup \Gamma$,
determine the source term $F(\cdot,t_0)$ in $\Omega_0$ by the observation data of 
$\{ (u,\nabla u, \nabla\mathrm{rot}\, u, H, \nabla H, \nabla\mathrm{rot}\, H)|_{\Gamma \times I}, \; 
u(\cdot,t_0)|_{\Omega}, \, H(\cdot, t_0)|_{\Omega} \}$.

By introducing the admissible set:
\begin{align*}
\mathcal{F}_{M,2} = 
\Big\{ 
& F\in H^2(0,T;H^2(\Omega)); 
F(\cdot,t_0)|_{\Gamma} = 0, 
\nabla F(\cdot,t_0)|_{\Gamma} = 0, \ 
\|\nabla F(\cdot,t_0)\|_{L^2(\partial\Omega)} \\
&+ \|F(\cdot,t_0)\|_{L^2(\partial\Omega)} \le M, 
\mbox{ \eqref{con:F2-2a} and \eqref{con:F2-2b} are satisfied}
\Big\}
\end{align*}
with arbitrarily given constant $M>0$ and the norm of the data: 
\begin{align*}
D_2 := \ 
&\|\nabla \mathrm{rot}\, u\|_{H^2(I;L^2(\Gamma))} 
+ \|\nabla \mathrm{rot}\, H\|_{H^2(I;L^2(\Gamma))} 
+ \|\nabla u\|_{H^3(I;L^2(\Gamma))} \\
&+ \|\nabla H\|_{H^3(I;L^2(\Gamma))} 
+ \|u\|_{H^2(I;L^2(\Gamma))} 
+ \|H\|_{H^2(I;L^2(\Gamma))} \\
&+ \|u(\cdot,t_0)\|_{H^4(\Omega)} 
+ \|H(\cdot,t_0)\|_{H^3(\Omega)},
\end{align*}
we have the following stability estimate. 
\begin{thm}
\label{thm:stab2}
Assume that $F\in \mathcal{F}_{M,2}$ and 
$(u,p,H)\in H^{2,3}(Q)\times H^{2,2}(Q)\times H^{2,3}(Q)$ 
is a solution to \eqref{sy:MHD} satisfying 
$\mathrm{rot}\, u,\mathrm{rot}\, H\in H^{2,3}(Q)$ and 
\begin{align*}
&\|u\|_{H^3(I;H^2(\Omega))} 
+ \|H\|_{H^3(I;H^2(\Omega))} 
+ \|\nabla \mathrm{rot}\, u\|_{H^2(I;L^2(\partial\Omega))} \\
&+ \|\nabla \mathrm{rot}\, H\|_{H^2(I;L^2(\partial\Omega))} 
+ \|\nabla u\|_{H^3(I;L^2(\partial\Omega))} 
+ \|\nabla H\|_{H^3(I;L^2(\partial\Omega))} \\
&+ \|u\|_{H^2(I;L^2(\partial\Omega))} 
+ \|H\|_{H^2(I;L^2(\partial\Omega))} \le M.
\end{align*} 
Then there exist constants $C > 0$ and $\theta\in (0,1)$ such that
\begin{align*}
\|F(\cdot,t_0)\| _{H^1(\Omega_0)} \le C (D_2^\theta + D_2).
\end{align*}
\end{thm}

\noindent {\bf Remark 3.}
In Theorem \ref{thm:stab2}, we can actually determine additionally 
the source term in the second equation of \eqref{sy:MHD} at the same time. 
That is, let the second equation be substituted by 
$$
\partial_t H - \kappa\Delta H + (A^{(3)}\cdot\nabla) H 
+ (H\cdot\nabla) A^{(4)} + L_2u = G(x,t).
$$
Then under suitable assumptions on $F$ and $G$, we can prove the stability estimate:
\begin{align*}
\|F(\cdot,t_0)\| _{H^1(\Omega_0)} + \|G(\cdot,t_0)\| _{H^1(\Omega_0)} 
\le C (D_2^\theta + D_2).
\end{align*}

\noindent {\bf Remark 4.}
In contrast to Theorem \ref{thm:stab1}, 
here we need not any data of $p$ since we assume the condition \eqref{con:F2-2b}. 
Thus, Theorem \ref{thm:stab2} indicates that one can recover 
a transverse field source, that is, a source whose divergence is zero,  
even if we do not measure the pressure $p$.  
Indeed, according to the Helmholtz decomposition:
\begin{equation}
\label{con:decomp}
F(x,t) = \nabla \widetilde{p}(x,t) + \widetilde F(x,t)
\end{equation}
with a real-valued function $\widetilde{p}$ and 
an $\mathbb{R}^3$-valued function $\widetilde F$ 
satisfying $\mbox{div}\, \widetilde F = 0$, 
we can rewrite the first equation in \eqref{sy:MHD} by
$$
\partial_t u - \nu\Delta u + (A^{(1)}\cdot\nabla) u 
+ (u\cdot\nabla) A^{(2)} + L_1H + \nabla (p-\widetilde{p}) 
= \widetilde F(x,t)                                            \quad \mbox{in} \ Q.
$$
Therefore, there are no hopes of determining the longitudinal field component $\nabla \widetilde{p}$ 
as long as we do not have the knowledge of $\nabla p$ and we can determine at most 
the information of the transverse field component $\widetilde F$. 

As we see from the above remark, 
we can identify only the transverse field component of the source 
provided that the pressure $p$ is not measured. 
In the next theorem, we consider the case that $F$ is not necessarily a transverse field, 
that is, we do not assume the condition \eqref{con:F2-2b}, 
and we determine the information of the transverse field component of the source. 

Note that the decomposition \eqref{con:decomp} is not unique. 
For example, we could replace $\nabla \widetilde p$ by $\nabla (\widetilde p-q)$ 
and $\widetilde F$ by $\widetilde F + \nabla q$ for any $q$ satisfying $\Delta q = 0$. 
Thus, it is not a suitable problem to recover $\widetilde F$ itself. 
Instead, we show the next theorem identifying 
$\mathrm{rot}\, \widetilde F = \mathrm{rot}\, F$ 
within the following admissible set:  
\begin{align*}
\mathcal{F}_{M,3} = 
\Big\{
F\in H^1(0,T;L^2(\Omega)); 
\mathrm{rot}\, F\in H^1(0,T;L^2(\Omega)), 
\mbox{ \eqref{con:F2-1} is satisfied}
\Big\}
\end{align*}
and under the norm of the data: 
\begin{align*}
D_3 := \ 
&\|\nabla \mathrm{rot}\, u\|_{H^1(I;L^2(\Gamma))} 
+ \|\nabla \mathrm{rot}\, H\|_{H^1(I;L^2(\Gamma))} 
+ \|\nabla u\|_{H^2(I;L^2(\Gamma))} \\
&+ \|\nabla H\|_{H^2(I;L^2(\Gamma))} 
+ \|u\|_{H^1(I;L^2(\Gamma))} 
+ \|H\|_{H^1(I;L^2(\Gamma))} 
+ \|u(\cdot,t_0)\|_{H^3(\Omega)}.
\end{align*}

\begin{thm}
\label{thm:stab3}
Assume that $F\in \mathcal{F}_{M,3}$ and 
$(u,p,H)\in H^{2,2}(Q)\times H^{2,1}(Q)\times H^{2,2}(Q)$ 
is a solution to \eqref{sy:MHD} satisfying 
$\mathrm{rot}\, u,\mathrm{rot}\, H\in H^{2,2}(Q)$ and 
\begin{align*}
&\|u\|_{H^2(I;H^2(\Omega))} 
+ \|H\|_{H^2(I;H^2(\Omega))} 
+ \|\nabla \mathrm{rot}\, u\|_{H^1(I;L^2(\partial\Omega))} \\
&+ \|\nabla \mathrm{rot}\, H\|_{H^1(I;L^2(\partial\Omega))} 
+ \|\nabla u\|_{H^2(I;L^2(\partial\Omega))} 
+ \|\nabla H\|_{H^2(I;L^2(\partial\Omega))} \\
&+ \|u\|_{H^1(I;L^2(\partial\Omega))} 
+ \|H\|_{H^1(I;L^2(\partial\Omega))} \le M.
\end{align*} 
Then there exist constants $C > 0$ and $\theta\in (0,1)$ such that
\begin{align*}
\|\mathrm{rot}\, F(\cdot,t_0)\| _{L^2(\Omega_0)} \le C (D_3^\theta + D_3).
\end{align*}
\end{thm}

\noindent {\bf Remark 5.}
Since the unknown source $F$ appears only in the first equation of \eqref{sy:MHD}, 
we can take advantage of the coupling between $u$ and $H$, 
and we manage to remove the measurement of $H(\cdot,t_0)$ in the data. 

\noindent {\bf Remark 6.}
Under some suitable boundary conditions on $\widetilde F(\cdot,t_0)$ (see e.g., \cite{Z06}), 
the decomposition \eqref{con:decomp} at $t=t_0$ can be uniquely determined. 
Hence one can recover the transverse field component $\widetilde F(\cdot,t_0)$ in $\Omega_0$ 
by using $\mathrm{rot}\, \widetilde F(\cdot,t_0) = \mathrm{rot}\, F(\cdot,t_0)$ in $\Omega_0$.

We end up this section with two examples which justify the conditions \eqref{con:F2-1} 
and \eqref{con:F2-2a}. 

\noindent {\bf Example 1.}
We assume that the temporal varying factor of the force $F$ is independent of the position $x$, 
that is, we let 
$$
F(x,t) = r(t) f(x),                                             \quad x\in \Omega,\ 0<t<T
$$
where $r\in C^2[0,T]$ is a given real-valued function and $f = (f_1,f_2,f_3)^{T}$,
$f_j\in C^1(\overline{\Omega})$, $j=1,2,3$. 
Then \eqref{con:F2-1} and \eqref{con:F2-2a} can be justified as long as 
$$
r(t_0) \not= 0.
$$
Actually, we have  
$$
|\partial_t^k \mathrm{rot}\, F(x,t)| 
= \left| \frac{d^k r}{dt^k}(t) \mathrm{rot}\, f(x)\right| 
\le \|r\|_{C^2[0,T]}|\mathrm{rot}\, f(x)| 
= \frac{\|r\|_{C^2[0,T]}}{|r(t_0)|} |\mathrm{rot}\, F(x,t_0)|,
$$
for all $(x,t)\in Q$ and $k=0,1,2$, which implies \eqref{con:F2-1}. 
By noting that 
$$
|\mathrm{rot}\, F(x,t_0)| \le C|\nabla F(x,t_0)|,                              \quad x\in \Omega,
$$
we immediately find \eqref{con:F2-2a}. 
Moreover, \eqref{con:F2-2b} holds true if we further assume $\mathrm{div}\, f = 0$. 

\noindent {\bf Example 2.}
We consider a more general case that the temporal varying factor of the force 
depends also on the position $x$, which can be modeled by a matrix $R(x,t)$: 
$$
F(x,t) = R(x,t)f(x),                                                   \quad x\in \Omega,\ 0<t<T,
$$
where $R = (r_{ij})_{1 \le i,j \le 3}$ is a given matrix-valued function and 
$f = (f_1,f_2,f_3)^T$ is a vector-valued function. 
Here we assume that each entry $r_{ij}$ and each component $f_j$ are smooth enough, 
for example,
$r_{ij}\in C^2([0,T];C^1(\overline{\Omega}))$ for $i,j=1,2,3$ and 
$f_j\in C^1(\overline{\Omega})$ for $j=1,2,3$.
Then \eqref{con:F2-2a} is satisfied if 
$$
|\mathrm{det}\, R(x,t_0)| > 0,                                 \quad x\in \overline{\Omega}. 
$$
In fact, by direct calculations, we have
\begin{align}
\label{ex:eq1}
|\partial_t^k \mathrm{rot}\, F(x,t)| 
\le C\left(|\nabla f(x)| + |f(x)|\right),                           \quad (x,t)\in Q,\ k=0,1,2.
\end{align}
On the other hand, we have
$$
\partial_\ell F(x,t) 
= \partial_\ell R(x,t) f(x) + R(x,t)\partial_\ell f(x),                              \quad (x,t)\in Q,\ \ell = 1,2,3.
$$
Since $|\mathrm{det}\, R(\cdot,t_0)| > 0$ on $\overline{\Omega}$, 
we find that $R(\cdot,t_0)$ is invertible and the norm of its inverse $R^{-1}(\cdot,t_0)$ 
is bounded from above. 
Thus, we rewrite the above equality at $t=t_0$ to obtain
$$
\partial_\ell f(x) 
= R^{-1}(x,t_0) \partial_\ell F(x,t_0) - R^{-1}(x,t_0) \partial_\ell R(x,t_0) f(x), 
                                                                                \quad x\in \Omega,\ \ell = 1,2,3,
$$
which implies
\begin{align}
\label{ex:eq2}
|\partial_\ell f(x)| \le C\left(|\partial_\ell F(x,t_0)| + |f(x)|\right), 
                                                                                \quad x\in \Omega,\ \ell = 1,2,3. 
\end{align}
Therefore, we obtain \eqref{con:F2-2a} by combining \eqref{ex:eq1} and \eqref{ex:eq2} with 
the following inequality
\begin{align*}
|f(x)| = |R^{-1}(x,t_0) F(x,t_0)| \le C|F(x,t_0)|,                        \quad x\in \Omega.
\end{align*}

We remark that in Example 2, we do not know a convenient sufficient condition 
for \eqref{con:F2-1} in terms of $R$ and $f$. 

The rest of this article is organized as follows. 
In Section 2, we show the key Carleman estimates which are stated in Theorems \ref{thm:CE1} 
and \ref{thm:CE2}. 
In Section 3, we prove Theorems \ref{thm:stab1}--\ref{thm:stab3} on the basis of 
these Carleman estimates. 
\section{Key Carleman estimates}
\label{sec-main}

In this section, we consider the following linearized MHD system 
\begin{equation}
\label{sy:ce}
\left\{
\begin{aligned}
&\ \partial_t u - \nu\Delta u + (A^{(1)}\cdot\nabla) u 
 + (u\cdot\nabla) A^{(2)} + L_1H + \nabla p = F            &&\quad \mbox{in} \ Q, \\
&\ \partial_t H - \kappa\Delta H + (A^{(3)}\cdot\nabla) H 
 + (H\cdot\nabla) A^{(4)} + L_2u = G                            &&\quad \mbox{in} \ Q, \\
&\ \mathrm{div}\, u = h                                                  &&\quad \mbox{in} \ Q.
\end{aligned}
\right.
\end{equation}
Here $L_1,L_2$ are defined in \eqref{def:L1L2} and the vector coefficients 
$A^{(j)}, C^{(k)}$, $j=1,2,3,4$, $k=1,2,3,4,5$ are smooth enough (e.g., $W^{3,\infty}(Q)$). 

The purpose of this section is to establish suitable Carleman estimates for \eqref{sy:ce}, 
which are the key points in proving our stability results of Theorems \ref{thm:stab1}--\ref{thm:stab3}. 
To this end, we first introduce the weight function for our Carleman estimates. 

Let $\Gamma \subset \partial\Omega$ be an arbitrarily fixed non-empty relatively open sub-boundary. 
We choose a bounded domain $\widetilde\Omega \subset \mathbb{R}^3\setminus \overline\Omega$ 
such that $\partial\widetilde\Omega \cap \partial\Omega = \overline\Gamma$. 
We construct a domain $\Omega_1 \supset \Omega$ by setting
$$
\Omega_1 = \Omega \cup \widetilde\Omega \cup \Gamma.
$$ 
Then 
\begin{equation*}
\overline\Gamma = \partial\Omega \cap \overline{\Omega_1}, \quad 
\partial\Omega_1 = 
(\partial\Omega \setminus \Gamma) \cup (\partial\widetilde\Omega \setminus \Gamma). 
\end{equation*}
With this new domain $\Omega_1$, we can choose a non-empty domain 
$\omega \subset \Omega_1$ such that 
$\overline\omega \subset \Omega_1 \setminus \overline\Omega = \widetilde\Omega$. 
Thus, referring to e.g., Fursikov and Imanuvilov \cite{FI96} or Imanuvilov \cite{I95}, 
we can construct a function $d\in C^2(\overline{\Omega_1})$ such that
\begin{equation}
\label{con:d}
d=0                 \quad \mbox{on $\partial\Omega_1$}, \quad 
d>0                 \quad \mbox{in $\Omega_1$}, \quad
|\nabla d| > 0   \quad \mbox{on } \overline{\Omega_1\setminus \omega}.
\end{equation}
In particular, $d>0$ in $\Omega$ and $|\nabla d|>0$ on $\overline\Omega$. 
For arbitrarily given $0<t_0<T$ and $0<\delta<\min\{t_0,T-t_0\}$, we denote
$$
I := (t_0-\delta,t_0+\delta), \quad 
Q := \Omega \times I, \quad 
Q_1 := \Omega_1 \times I. 
$$
Moreover, we fix a large constant $\lambda>0$, 
and for a constant $\beta>0$ we set 
\begin{equation}
\label{con:psi-phi}
\psi(x,t) = d(x) - \beta (t-t_0)^2, \quad 
\varphi(x,t) = e^{\lambda\psi(x,t)},                                      \quad (x,t)\in \overline{Q_1}.               
\end{equation}
Now we establish our Carleman estimates in the following two cases. 

\subsection{Carleman estimate with pressure term $p$}

Here we state the first Carleman estimate in which we include the estimate of $p$. 
\begin{thm}
\label{thm:CE1}
Let $F\in L^2(Q)$. 
Then there exist constants $s_0>0$ and $C>0$ such that 
\begin{align*}
&\int_Q \Biggl( 
\frac{1}{s^2}\left( \vert \partial_t u\vert^2 + \vert \Delta u\vert^2\right)
+ \frac{1}{s}\left(\vert \partial_t H\vert^2 + \vert \Delta H\vert^2\right) 
+ \vert \nabla u\vert^2 
+ s\vert \nabla H\vert^2 \\
&+ s^2\vert u\vert^2 
+ s^3\vert H\vert^2
+ \frac{1}{s}\vert \nabla p\vert^2 + s\vert p\vert^2
\Biggr) e^{2s\varphi} dxdt \\
\le &\; C\int_Q 
\left(\vert F\vert^2 + \vert G\vert^2 + \vert \nabla_{x,t}h\vert^2 \right) 
e^{2s\varphi} dxdt \\
&+ Cs^3\int_{\partial\Omega \times I} 
(\vert u\vert^2 + \vert H\vert^2 + \vert \nabla_{x,t}u\vert^2 + \vert \nabla_{x,t}H\vert^2) 
e^{2s\varphi} dSdt \\
&+ Cs^3\int_\Omega 
(\vert u\vert^2 + \vert H\vert^2 + \vert \nabla u\vert^2 + \vert \nabla H\vert^2)
e^{2s\varphi} dx\Big|_{t=t_0\pm \delta}
\end{align*}
for all $s\ge s_0$ and 
\begin{equation*}
u,H\in H^1(I;H^1(\Omega)) \cap L^2(I;H^2(\Omega)), \quad 
p\in L^2(I;H_0^1(\Omega)) 
\end{equation*}
satisfying \eqref{sy:ce}. 
\end{thm}
Henceforth $C>0$ denotes a generic constant which depends on 
$\lambda, t_0, \delta, T, \beta, \Omega$ and the coefficients, 
but is independent of the parameter $s$. 
Such uniformity in $s$ is used for deriving the stability results for the inverse problems. 

\begin{proof}[Proof of Theorem \ref{thm:CE1}]
We divide the proof into three steps. 

\noindent{\bf First step:}

At the beginning of the proof, we assume that $p\in L^2(I;C_0^\infty(\Omega))$. 
Later in the third step, we will take approximations and treat the case of 
$p\in L^2(I;H_0^1(\Omega))$. 
In the first step, we prove the following lemma. 
\begin{lem}
\label{lem-p}
Let $p\in L^2(I;C_0^\infty(\Omega))$ satisfy
\begin{equation}
\label{eq:p}
\Delta p = f_0 + \sum_{j=1}^3 \partial_j f_j                                 \quad \mbox{in }Q
\end{equation}
Then there exist constants $\hat s>0$ and $C>0$ such that 
$$
\int_Q (|\nabla p|^2 + s^2|p|^2) e^{2s\varphi} dxdt 
\le C\int_Q \left(\frac{1}{s}|f_0|^2 + s\sum_{j=1}^3 |f_j|^2 \right) e^{2s\varphi} dxdt
$$
for all $s\ge \hat s$. 
\end{lem}
\begin{proof}
By noting that $p\in L^2(I;C_0^\infty(\Omega))$, 
we can take the zero extension of $p$ from $\Omega$ to $\Omega_1$, 
which is denoted by the same letter.  
Then we find that $p(\cdot,t) \in C_0^\infty(\Omega_1)$ 
for $t\in I$ and \eqref{eq:p} holds true in $Q_1$ as long as we take also the zero extensions of 
$f_j$, $j=0,1,2,3$. 

Now we apply the $H^{-1}$ Carleman estimate for the elliptic equation 
(Theorem A.1 in Imanuvilov and Puel \cite{IP03}) to \eqref{eq:p} in $\Omega_1$ 
and we obtain 
\begin{align*}
&\int_{\Omega_1} (|\nabla p(x,t)|^2 + s^2|p(x,t)|^2) e^{2se^{\lambda d(x)}} dx \\
\le &\; C\int_{\Omega_1} \left(\frac{1}{s}|f_0(x,t)|^2 
 + s\sum_{j=1}^3|f_j(x,t)|^2 \right) e^{2se^{\lambda d(x)}} dx 
+ C\int_\omega s^2|p(x,t)|^2 e^{2se^{\lambda d(x)}}dx
\end{align*}
for all $s\ge s_1$, $t\in I$ and $\overline\omega \subset \Omega_1$. 
Recalling the construction of the domain $\Omega_1$, 
we can choose a non-empty domain $\omega$ such that 
$\overline\omega \subset \Omega_1\setminus \overline\Omega$. 
Since we take the zero extensions for $f_j$, $p$ outside of $\Omega$, 
the integral over $\omega$ vanishes and we have
\begin{align}
\nonumber
&\int_{\Omega} (|\nabla p(x,t)|^2 + s^2|p(x,t)|^2) e^{2se^{\lambda d(x)}} dx \\
\label{ce:eq1}
\le &\; C\int_{\Omega} \left(\frac{1}{s}|f_0(x,t)|^2 
 + s\sum_{j=1}^3|f_j(x,t)|^2 \right) e^{2se^{\lambda d(x)}} dx
\end{align}
for all $s\ge s_1$ and each $t\in I$. 
Let $\hat s := s_1 e^{\lambda\beta \delta^2}$. 
Then $s\ge \hat s$ implies 
$$
se^{-\lambda\beta (t-t_0)^2} 
\ge se^{-\lambda\beta \delta^2} 
\ge \hat se^{-\lambda\beta \delta^2} 
= s_1,
$$
for $t\in I$, and thus, replacing $s$ by $se^{-\lambda\beta (t-t_0)^2}$ in \eqref{ce:eq1} yields
\begin{align*}
&\int_{\Omega} (|\nabla p(x,t)|^2 + s^2e^{-2\lambda\beta (t-t_0)^2}|p(x,t)|^2) e^{2s\varphi(x,t)} dx \\
\le &\; C\int_{\Omega} \left(\frac{1}{s}e^{\lambda\beta (t-t_0)^2}|f_0(x,t)|^2 
 + se^{-\lambda\beta (t-t_0)^2}\sum_{j=1}^3|f_j(x,t)|^2 \right) e^{2s\varphi(x,t)} dx
\end{align*}
for all $s\ge \hat s$ and each $t\in I$. 
Here we used $\varphi(x,t) = e^{\lambda d(x)} e^{-\lambda\beta (t-t_0)^2}$, $(x,t)\in Q$. 
Therefore, we complete the proof of Lemma \ref{lem-p} by integrating the above estimate over $t\in I$. 
\end{proof}

\noindent{\bf Second step:}

In this step, we apply the estimate in the first step and the Carleman estimate for parabolic system 
to prove the theorem in the case where $p\in L^2(I;C_0^\infty(\Omega))$. 

We derive the following elliptic equation with respect to $p$ by taking the divergence 
on both sides of the first equation in \eqref{sy:ce} and using the third equation there:  
\begin{align}
\label{ce:eq2}
\Delta p 
= \mathrm{div}\, F - \partial_t h + \nu \mathrm{div}(\nabla h) 
 - \mathrm{div}((A^{(1)}\cdot\nabla)u) - \mathrm{div}((u\cdot\nabla)A^{(2)}) 
 - \mathrm{div}(L_1 H).
\end{align}
In particular, we rewrite several terms on the right-hand side as follows
\begin{align*}
&\nu \mathrm{div}(\nabla h)
=\nu \sum_{j=1}^3 \partial_j (\partial_j h), \\
&\mathrm{div}((A^{(1)}\cdot\nabla)u) 
= \sum_{j,k=1}^3 \partial_k (A^{(1)}_j \partial_j u_k) 
= \sum_{j=1}^3 A^{(1)}_j \partial_j \left(\sum_{k=1}^3 \partial_k u_k\right) 
 + \sum_{j,k=1}^3 (\partial_k A^{(1)}_j) \partial_j u_k\\
&\hspace{2.4cm} = A^{(1)}\cdot\nabla (\mathrm{div}\, u) 
 + \sum_{j,k=1}^3\left( \partial_j((\partial_k A^{(1)}_j)u_k) - (\partial_j\partial_k A^{(1)}_j)u_k \right), \\
&\mathrm{div}((u\cdot\nabla)A^{(2)}) + \mathrm{div}(L_1 H)
= \sum_{j=1}^3 \partial_j((u\cdot\nabla)A^{(2)}_j + [L_1 H]_j) 
\end{align*}
where $[v]_j$ denotes the $j$-th component of a vector $v$, $j=1,2,3$. 
Besides, we have
\begin{align*}
\mathrm{div}\, F = \sum_{j=1}^3 \partial_j F_j.
\end{align*}
Therefore, we are able to represent the right-hand side of \eqref{ce:eq2} in the form of 
$f_0 + \sum_{j=1}^3 \partial_j f_j$ with
\begin{align*}
&f_0 
= - \partial_t h - (A^{(1)}\cdot\nabla) h 
 + \sum_{k,\ell=1}^3 (\partial_\ell\partial_k A^{(1)}_\ell)u_k,\\
&f_j 
= F_j + \nu \partial_j h 
- \sum_{k=1}^3 (\partial_k A^{(1)}_j)u_k 
- (u\cdot\nabla)A^{(2)}_j - [L_1 H]_j,                                                           \quad j=1,2,3. 
\end{align*}
Thus, by noting 
\begin{align*}
&|f_0| \le C \left( |\partial_t h| + |\nabla h| + |u|\right),\\
&|f_j| \le C \left(|F| + |\nabla h| + |u| + |H| + |\nabla H|\right),                      \quad j=1,2,3,
\end{align*}
we apply Lemma \ref{lem-p} to obtain
\begin{align}
\nonumber
\int_Q \left(\frac{1}{s}|\nabla p|^2 + s|p|^2\right) e^{2s\varphi} dxdt 
\le &\; C\int_Q \left(|F|^2 + |\nabla_{x,t} h|^2\right) e^{2s\varphi} dxdt\\
\label{ce:eq3}
&+ C\int_Q (|u|^2 + |H|^2 + |\nabla H|^2)e^{2s\varphi} dxdt
\end{align}
for all $s\ge s_2:= \max\{1,\hat s\}$. 

On the other hand, we can regard $\nabla p$ in the first equation of \eqref{sy:ce} as 
a non-homogeneous term, and then the first and second equations of \eqref{sy:ce} 
become a weakly coupling parabolic system with respect to $u$ and $H$. 
Hence we can apply the Carleman estimate for the parabolic system 
(e.g., Yamamoto \cite[Theorem 3.2]{Y09}) and we obtain
\begin{align}
\nonumber
&\int_Q \left(\frac{1}{s^2}(|\partial_t u|^2 + |\Delta u|^2) 
+ |\nabla u|^2 + s^2|u|^2 \right)e^{2s\varphi} dxdt\\
\nonumber
\le &\; C\int_Q \frac{1}{s}|F|^2e^{2s\varphi} dxdt 
+ C\int_Q \frac{1}{s}(|\nabla p|^2 + |\nabla H|^2 + |H|^2)e^{2s\varphi} dxdt\\
\label{ce:eq4}
&+ Cs^2\int_{\partial\Omega \times I} (|\nabla_{x,t} u|^2 + |u|^2)e^{2s\varphi} dSdt
+ Cs^2\int_\Omega (|\nabla u|^2 + |u|^2)e^{2s\varphi} dx \Big|_{t=t_0\pm \delta} 
\end{align}
and 
\begin{align}
\nonumber
&\int_Q \left(\frac{1}{s}(|\partial_t H|^2 + |\Delta H|^2) 
+ s|\nabla H|^2 + s^3|H|^2 \right)e^{2s\varphi} dxdt\\
\nonumber
\le &\; C\int_Q |G|^2e^{2s\varphi} dxdt 
+ C\int_Q (|\nabla u|^2 + |u|^2)e^{2s\varphi} dxdt\\
\label{ce:eq5}
&+ Cs^3\int_{\partial\Omega \times I} (|\nabla_{x,t} H|^2 + |H|^2)e^{2s\varphi} dSdt
+ Cs^3\int_\Omega (|\nabla H|^2 + |H|^2)e^{2s\varphi} dx \Big|_{t=t_0\pm \delta} 
\end{align}
for all $s\ge s_3$. 
Substituting \eqref{ce:eq5} into the last term on the right-hand side of \eqref{ce:eq3} 
and then substituting the derived inequality into the second term 
on the right-hand side of \eqref{ce:eq4}, we obtain
\begin{align}
\nonumber
&\int_Q \left(\frac{1}{s^2}(|\partial_t u|^2 + |\Delta u|^2) 
+ |\nabla u|^2 + s^2|u|^2 \right)e^{2s\varphi} dxdt\\
\nonumber
\le &\; C\int_Q \left(|F|^2 + \frac{1}{s}|G|^2 + |\nabla_{x,t} h|^2\right)e^{2s\varphi} dxdt \\
\nonumber
&+ C\int_Q \left(\frac{1}{s}|\nabla u|^2 + |u|^2 
 + \frac{1}{s}|\nabla H|^2 + \frac{1}{s}|H|^2 \right)e^{2s\varphi} dxdt\\
\nonumber
&+ Cs^2\int_{\partial\Omega \times I} (|\nabla_{x,t} u|^2 + |u|^2 
 + |\nabla_{x,t} H|^2 + |H|^2)e^{2s\varphi} dSdt\\
\label{ce:eq6}
&+ Cs^2\int_\Omega (|\nabla u|^2 + |u|^2 
 + |\nabla H|^2 + |H|^2)e^{2s\varphi} dx \Big|_{t=t_0\pm \delta} 
\end{align}
for all $s\ge s_4 := \max\{s_2,s_3\}$. 
Therefore, we combine the inequalities \eqref{ce:eq3}, \eqref{ce:eq5} and \eqref{ce:eq6}, 
by choosing $s>0$ sufficiently large, we can absorb the lower-order terms 
$\int_Q (|\nabla u|^2 + s|u|^2 + |\nabla H|^2 + |H|^2)e^{2s\varphi} dxdt$ 
into the left-hand side and this completes the proof in the case of $p\in L^2(I;C_0^\infty(\Omega))$.

\noindent{\bf Third step:}

In this step, we consider the case of $p\in L^2(I; H_0^1(\Omega))$ by a density argument. 
Recall the density $\overline{C_0^\infty(\Omega)} = H_0^1(\Omega)$ by the norm of $H^1(\Omega)$ 
(e.g., Adams and Fournier\cite{AF03}). 
This yields the existence of an approximating sequence 
$p_n\in L^2(I;C_0^\infty(\Omega))$, $n\in \mathbb{N}$ such that 
\begin{align}
\label{ce:eq7}
\lim_{n\to \infty} \|p_n-p\|_{L^2(I;H^1(\Omega))} = 0.
\end{align}
Then we rewrite the first equation of \eqref{sy:ce} in terms of $p_n$: 
$$
\partial_t u - \nu \Delta u + (A\cdot\nabla)u + (u\cdot\nabla)B + L_1 H + \nabla p_n 
= F + \nabla(p_n-p).
$$
Therefore, by applying the Carleman estimate in the above steps, we have
\begin{align*}
&\int_Q \Biggl( 
\frac{1}{s^2}\left( \vert \partial_t u\vert^2 + \vert \Delta u\vert^2 \right)
+ \frac{1}{s}\left( \vert \partial_t H\vert^2 + \vert \Delta H\vert^2\right) 
+  \vert \nabla u\vert^2 
+ s \vert \nabla H\vert^2 \\
&+ s^2\vert u\vert^2 
+ s^3\vert H\vert^2
+ \frac{1}{s}\vert \nabla p\vert^2 
+ s\vert p\vert^2
\Biggr) e^{2s\varphi} dxdt \\
\le &\; C\int_Q 
\left(\vert F + \nabla(p_n-p)\vert^2 + \vert G\vert^2 + \vert \nabla_{x,t}h\vert \right) 
e^{2s\varphi} dxdt \\
&+ Cs^3\int_{\partial\Omega \times I} 
(\vert u\vert^2 + \vert H\vert^2 + \vert \nabla_{x,t}u\vert^2 + \vert \nabla_{x,t}H\vert^2) 
e^{2s\varphi} dSdt \\
&+ Cs^3\int_\Omega 
(\vert u\vert^2 + \vert H\vert^2 + \vert \nabla u\vert^2 + \vert \nabla H\vert^2)
e^{2s\varphi} dx\Big|_{t=t_0\pm \delta}
\end{align*}
for all $s\ge s_0$. Since 
$$
|F + \nabla(p_n-p)|^2 
\le C\left(|F|^2 + |\nabla(p_n-p)|^2\right),
$$
by letting $n\to \infty$ and using \eqref{ce:eq7}, we complete the proof of Theorem \ref{thm:CE1}.
\end{proof}
\subsection{Carleman estimate without pressure term $p$}

In this subsection, we propose another Carleman estimate for the linearized MHD system 
\eqref{sy:ce} in which we do not estimate $p$.
Such estimate is useful in proving the stability results without the information of the pressure. 

Although we can argue with less regularity which is similar to the previous subsection, 
here we assume higher regularity for the non-homogeneous term and the coefficients 
in order to clarify the essential idea. 

Recall that $\Gamma \subset \partial\Omega$ is an arbitrarily given 
non-empty relatively open sub-boundary, the functions $d$ and $\varphi$ 
are given by \eqref{con:d} and \eqref{con:psi-phi} respectively. 
Now we state the second Carleman estimate. 
\begin{thm}
\label{thm:CE2}
Let $h=0$, $F,G\in L^2(Q)$ and 
$\mathrm{rot}\, F, \mathrm{rot}\, G \in L^2(Q)$. 
Then there exist constants $C>0$ and $\widehat{s_0}>0$ such that
\begin{align*}
&\int_Q \bigg(
\frac{1}{s}\left(|\partial_t \mathrm{rot}\, u|^2 + |\Delta \mathrm{rot}\, u|^2
 + |\partial_t \mathrm{rot}\, H|^2 + |\Delta \mathrm{rot}\, H|^2\right) \\
&+ s\left(|\nabla \mathrm{rot}\, u|^2 + |\Delta u|^2 + |\nabla \mathrm{rot}\, H|^2 
 + |\Delta H|^2 \right) + s^2\left(|\nabla u|^2 + |\nabla H|^2\right) \\
&+ s^3\left(|\mathrm{rot}\, u|^2 + |\mathrm{rot}\, H|^2 \right) 
 + s^4\left(|u|^2 + |H|^2\right) 
\bigg)e^{2s\varphi} dxdt\\
\le &\; C\int_Q \left(|\mathrm{rot}\, F|^2 + |\mathrm{rot}\, G|^2\right)e^{2s\varphi} dxdt \\
&+ Cs^4\int_{\partial\Omega \times I} \left(|\nabla_{x,t} \mathrm{rot}\, u|^2 
 + |\nabla u|^2 + |\nabla_{x,t} \mathrm{rot}\, H|^2 + |\nabla H|^2 
 + |u|^2 + |H|^2 \right)e^{2s\varphi} dSdt\\
&+ Cs^3\int_\Omega \left(|\nabla \mathrm{rot}\, u|^2 + |\mathrm{rot}\, u|^2 
 + |\nabla \mathrm{rot}\, H|^2 + |\mathrm{rot}\, H|^2 \right)e^{2s\varphi} dx \Big|_{t=t_0\pm \delta} 
\end{align*}
for all $s\ge \widehat{s_0}$ and 
\begin{equation*}
u,H\in H^1(I;H^2(\Omega)), \quad 
\mathrm{rot}\, u,\mathrm{rot}\, H\in H^1(I;H^2(\Omega)), \quad 
p\in L^2(I;H^2(\Omega)) 
\end{equation*}
satisfying \eqref{sy:ce} and $\mathrm{div}\, H = 0$.
\end{thm}
\begin{proof}[Proof of Theorem \ref{thm:CE2}]
The key to the proof is the reduction of \eqref{sy:ce} to a parabolic-elliptic system. 
Then we can apply standard Carleman estimates for parabolic equations and elliptic equations 
to establish our desired Carleman estimate. 

To start with, we denote $y := \mathrm{rot}\, u$, $z := \mathrm{rot}\, H$. 
By noting 
$$
\Delta u = -\mathrm{rot}(\mathrm{rot}\, u) + \nabla (\mathrm{div}\, u) \quad 
\mbox{and }\quad
\mathrm{rot}(\nabla p) = 0,
$$
the system \eqref{sy:ce} and $\mathrm{div}\, H = 0$ yield
\begin{equation}
\label{ce2:sy1}
\left\{
\begin{aligned}
&\ \partial_t y - \nu\Delta y + (A^{(1)}\!\cdot\!\nabla) y 
=  \mathrm{rot}\, F -\! (C^{(1)}\!\cdot\nabla) z 
 -\! \sum_{k=1}^3 \nabla A^{(1)}_k \times \partial_k u 
 -\! (u\!\cdot\!\nabla)\mathrm{rot}\, A^{(2)} \\
 &\hspace{0.7cm} -\! \sum_{k=1}^3 \nabla u_k \times \partial_k A^{(2)} 
 -\! \sum_{k=1}^3 \nabla C^{(1)}_k \times \partial_k H 
 -\! (H\!\cdot\!\nabla) \mathrm{rot}\, C^{(2)} 
 -\! \sum_{k=1}^3 \nabla H_k \times \partial_k C^{(2)},\\
&\ \partial_t z - \kappa\Delta z + (A^{(3)}\!\cdot\!\nabla) z 
= \mathrm{rot}\, G -\! (C^{(5)}\!\cdot\nabla) y
 -\! \sum_{k=1}^3 \nabla A^{(3)}_k \times \partial_k H 
 -\! (H\!\cdot\!\nabla)\mathrm{rot}\, A^{(4)} \\
 &\hspace{0.7cm} -\! \sum_{k=1}^3 \nabla H_k \times \partial_k A^{(4)} 
 -\! \sum_{k=1}^3 \nabla C^{(5)}_k \times \partial_k u 
 -\! (u\!\cdot\!\nabla) \mathrm{rot}\, C^{(4)} 
 -\! \sum_{k=1}^3 \nabla u_k \times \partial_k C^{(4)},\\
&\ \Delta u = -\mathrm{rot}\, y, \\
&\ \Delta H = -\mathrm{rot}\, z.
\end{aligned}
\right.
\end{equation}
Next we employ the Carleman estimate for the parabolic system 
(e.g., Yamamoto \cite[Theorem 3.2]{Y09}) 
to the first and the second equations of \eqref{ce2:sy1} and we obtain
\begin{align}
\nonumber
&\int_Q \left(\frac{1}{s}(|\partial_t y|^2 + |\Delta y|^2) 
 + s|\nabla y|^2 + s^3|y|^2 \right)e^{2s\varphi} dxdt\\
\nonumber
\le &\; C\int_Q |\mathrm{rot}\, F|^2e^{2s\varphi} dxdt 
 + C\int_Q (|\nabla u|^2 + |u|^2 + |\nabla H|^2 + |H|^2 + |\nabla z|^2)e^{2s\varphi} dxdt\\
\label{ce2:eq1}
&+ Cs^3\int_{\partial\Omega \times I} (|\nabla_{x,t} y|^2 + |y|^2)e^{2s\varphi} dSdt
 + Cs^3\int_\Omega (|\nabla y|^2 + |y|^2)e^{2s\varphi} dx \Big|_{t=t_0\pm \delta} 
\end{align}
and 
\begin{align}
\nonumber
&\int_Q \left(\frac{1}{s}(|\partial_t z|^2 + |\Delta z|^2) 
 + s|\nabla z|^2 + s^3|z|^2 \right)e^{2s\varphi} dxdt\\
\nonumber
\le &\; C\int_Q |\mathrm{rot}\, G|^2e^{2s\varphi} dxdt 
 + C\int_Q (|\nabla u|^2 + |u|^2 + |\nabla H|^2 + |H|^2 + |\nabla y|^2)e^{2s\varphi} dxdt\\
\label{ce2:eq2}
&+ Cs^3\int_{\partial\Omega \times I} (|\nabla_{x,t} z|^2 + |z|^2)e^{2s\varphi} dSdt
 + Cs^3\int_\Omega (|\nabla z|^2 + |z|^2)e^{2s\varphi} dx \Big|_{t=t_0\pm \delta} 
\end{align}
for all $s\ge \widehat{s_1}$. 
On the other hand, we employ the following Carleman estimate for the elliptic equation. 
\begin{lem}
\label{lem:ce-ellip}
Let $g\in L^2(Q)$ and $w\in L^2(I;H^2(\Omega))$ satisfy
$$
\Delta w(x,t) = g(x,t),                               \quad (x,t)\in Q.
$$
Then there exist constants $\widehat{s_2}>0$ and $C>0$ such that
\begin{align*}
&\int_Q \left(|\Delta w|^2 + s|\nabla w|^2 + s^3|w|^2\right)e^{2s\varphi} dxdt \\
\le &\; C\int_Q |g|^2 e^{2s\varphi} dxdt 
 + Cs^3\int_{\partial\Omega \times I} (|\nabla w|^2 + |w|^2)e^{2s\varphi} dSdt 
\end{align*}
for all $s\ge \widehat{s_2}$. 
\end{lem}
The proof is based on the standard elliptic Carleman estimate which is well-known. 
For example, it can be directly proved by following similar steps as 
the proof of \cite[Theorem 3.1]{Y09}, and keeping the boundary integrals. 
 
By applying Lemma \ref{lem:ce-ellip} to the third and the fourth equations of \eqref{ce2:sy1} 
and multiplying both sides by $s$, we have
\begin{align}
\nonumber
&\int_Q \left( s|\Delta u|^2 + s^2|\nabla u|^2 + s^4|u|^2\right)e^{2s\varphi} dxdt \\
\label{ce2:eq3}
\le &\; C\int_Q s|\mathrm{rot}\, y|^2 e^{2s\varphi} dxdt 
 + Cs^4\int_{\partial\Omega \times I} (|\nabla u|^2 + |u|^2)e^{2s\varphi} dSdt
\end{align}
and 
\begin{align}
\nonumber
&\int_Q \left( s|\Delta H|^2 + s^2|\nabla H|^2 + s^4|H|^2\right)e^{2s\varphi} dxdt \\
\label{ce2:eq4}
\le &\; C\int_Q s|\mathrm{rot}\, z|^2 e^{2s\varphi} dxdt 
+ Cs^4\int_{\partial\Omega \times I} (|\nabla H|^2 + |H|^2)e^{2s\varphi} dSdt 
\end{align}
for all $s\ge \widehat{s_2}$. 
Since 
$$
|\mathrm{rot}\, y|^2 \le C|\nabla y|^2,\quad 
|\mathrm{rot}\, z|^2 \le C|\nabla z|^2, 
$$
we sum up the estimates \eqref{ce2:eq1}--\eqref{ce2:eq4} 
and insert \eqref{ce2:eq1}--\eqref{ce2:eq2} into the right-hand side, 
and then we can find a sufficiently large constant 
$\widehat{s_0}\ge \max\{\widehat{s_1},\widehat{s_2}\}$ 
such that we absorb the term 
$C\int_Q (|\nabla u|^2 + |u|^2 + |\nabla H|^2 + |H|^2 + |\nabla y|^2 + |\nabla z|^2)e^{2s\varphi} dxdt$ 
on the right-hand side into the left-hand side. 
This completes the proof of Theorem \ref{thm:CE2}. 
\end{proof}
\section{Proofs of Theorems \ref{thm:stab1}--\ref{thm:stab3}}
\subsection{Proof of Theorem \ref{thm:stab1}} 
\label{subsec:proof1}
We divide the proof into four steps. 

\noindent {\bf First step: Choice of a suitable cut-off function $\chi$}

In Theorem \ref{thm:stab1}, we do not have $p=0$ on $\partial\Omega$, 
which prevents us from applying the Carleman estimate (Theorem \ref{thm:CE1}) 
established in the former section directly. 
Thus, here we need a cut-off argument. 

Recall that $d=d(x)$ satisfies \eqref{con:d}. 
Since $d(x)>0$ for $x\in \Omega_0 \subset \Omega_1$, we can choose constants 
$0<\varepsilon_1<\varepsilon_2<\varepsilon_3$ such that 
$\Omega_0 \subset \{x\in \Omega; d(x)>\varepsilon_3\}$. 

We choose a cut-off function as follows: 
$\chi\in C^{\infty}(\mathbb{R}^3)$ satisfying
$0 \le \chi\le 1$ and 
\begin{equation}
\label{con:chi}
\chi(x) = 
\left\{
\begin{aligned}
& 1\quad \mbox{if } d(x)>\varepsilon_2,\                        x\in \overline\Omega, \\
& 0\quad \mbox{if } d(x)<\varepsilon_1,\                        x\in \overline\Omega.
\end{aligned}
\right.
\end{equation}
By \eqref{con:d}, we have $d=0$ on $\partial\Omega \setminus \Gamma$. 
From \eqref{con:chi} and the continuity of $d$, 
we easily find that $\chi$ and its derivatives vanish on 
$\partial\Omega \setminus \Gamma$. 
On the other hand, by the Sobolev extension theorem, 
we can find $p_0\in H^1(I;H^1(\Omega))$ such that
\begin{equation}
\label{stab:eq1}
p_0 = p \quad                          \mbox{on }\Gamma \times I \quad 
\mbox{and}\quad 
\|p_0\|_{H^1(I;H^1(\Omega))} \le C\|p\|_{H^1(I;H^{\frac12}(\Gamma))}.
\end{equation}
Therefore, we set 
$$
q= (p-p_0),\quad 
q_\chi = \chi q,
$$ 
so that $q_\chi=0$ on $\partial\Omega \times I$ and 
$q_\chi\in H^1(I;H_0^1(\Omega))$. 
Then by letting 
$$
u_\chi = \chi u,\quad 
H_\chi = \chi H,
$$
we can rewrite the linearized MHD system \eqref{sy:MHD} by
\begin{equation}
\label{stab:sy1}
\left\{
\begin{aligned}
&\partial_t u_\chi - \nu\Delta u_\chi + (A^{(1)}\cdot\nabla) u_\chi 
 + (u_\chi\cdot\nabla) A^{(2)} + L_1 H_\chi + \nabla q_\chi 
= \chi F - \chi \nabla p_0 + J_1 \\
&\partial_t H_\chi - \kappa\Delta H_\chi + (A^{(3)}\cdot\nabla) H_\chi 
 + (H_\chi\cdot\nabla) A^{(4)} + L_2 u_\chi = J_2 \\
&\mathrm{div}\, u_\chi = \nabla\chi\!\cdot\! u
\end{aligned}
\right.
\end{equation}
Here $J_1$, $J_2$ include the terms of $p, p_0$ and 
at most first spatial derivatives of $u$, $H$ multiplied by 
the derivatives of the cut-off function $\chi$, in other words, we have 
\begin{align*}
|J_k| 
\le C\left(|\nabla \chi| + |\Delta \chi|\right)
 \left(|p| + |p_0| + |u| + |H| + |\nabla u| + |\nabla H| \right),
\end{align*}  
for $k=1,2$. 
Moreover, in order to prove the stability result, 
we need to differentiate the system once with respect to $t$. 
Henceforth, we denote
$$
u^{(1)} = \partial_t u, \ 
H^{(1)} = \partial_t H, \ 
q^{(1)} = \partial_t q, \quad 
u^{(1)}_\chi = \chi \partial_t u, 
\ H^{(1)}_\chi = \chi \partial_t H, 
\ q^{(1)}_\chi = \chi \partial_t q.
$$
Hence we can derive the system with respect to $u^{(1)}_\chi, H^{(1)}_\chi, q^{(1)}_\chi$ 
by \eqref{stab:sy1}: 
\begin{equation}
\label{stab:sy2}
\left\{
\begin{aligned}
&\partial_t u^{(1)}_\chi - \nu\Delta u^{(1)}_\chi + (A^{(1)}\cdot\nabla) u^{(1)}_\chi 
 + (u^{(1)}_\chi\cdot\nabla) A^{(2)} + L_1 H^{(1)}_\chi + \nabla q^{(1)}_\chi \\
&\quad = \chi (\partial_t F) - \chi \nabla (\partial_t p_0) + J_3 \\
&\partial_t H^{(1)}_\chi - \kappa\Delta H^{(1)}_\chi + (A^{(3)}\cdot\nabla) H^{(1)}_\chi 
 + (H^{(1)}_\chi\cdot\nabla) A^{(4)} + L_2 u^{(1)}_\chi = J_4 \\
&\mathrm{div}\, u^{(1)}_\chi = \nabla\chi\!\cdot\! \partial_t u
\end{aligned}
\right.
\end{equation}
Here $J_3$, $J_4$ include the terms of $p, p_0, \partial_t p, \partial_t p_0$ and 
at most first spatial derivatives of $u$, $H$, $\partial_t u$, $\partial_t H$ multiplied by 
the derivatives of the cut-off function $\chi$, in other words, we have 
\begin{align*}
|J_k| \le 
&\; C\left(|\nabla \chi| + |\Delta \chi|\right)
 \left(|p| + |p_0| + |u| + |H| + |\nabla u| + |\nabla H| \right)\\
&+ C\left(|\nabla \chi| + |\Delta \chi|\right)
 \left(|\partial_t p| + |\partial_t p_0| + |\partial_t u| + |\partial_t H| 
 + |\nabla \partial_t u| + |\nabla \partial_t H| \right),
\end{align*}  
for $k=3,4$.

\noindent{\bf Second step: Application of the Carleman estimate} 

Next we employ Theorem \ref{thm:CE1} to \eqref{stab:sy1}. 
In particular, we obtain
\begin{align}
\nonumber
&\int_Q \Biggl( 
\frac{1}{s^2}|\partial_t u_\chi|^2 + |\nabla u_\chi|^2 
+ s|\nabla H_\chi|^2 + s^2|u_\chi|^2 + s^3|H_\chi|^2
\Biggr) e^{2s\varphi} dxdt \\
\nonumber
\le &\; C\int_Q \left(|\chi F - \chi \nabla p_0 + J_1|^2 
 + |J_2|^2 + |\nabla_{x,t}(\nabla \chi \cdot u)|^2 \right) e^{2s\varphi} dxdt \\
\nonumber
&+ Cs^3\int_{\partial\Omega \times I} 
(|u_\chi|^2 + |H_\chi|^2 + |\nabla_{x,t}u_\chi|^2 + |\nabla_{x,t}H_\chi|^2) 
e^{2s\varphi} dSdt \\
\label{stab:eq2}
&+ Cs^3\int_\Omega 
(|u_\chi|^2 + |H_\chi|^2 + |\nabla u_\chi|^2 + |\nabla H_\chi|^2)
e^{2s\varphi} dx\Big|_{t=t_0\pm \delta}
\end{align}
for all sufficiently large $s\ge 1$. 
We estimate the right-hand side of the inequality \eqref{stab:eq2} in more details. 

For the first integral, by \eqref{con:F1} and using the triangle inequality, we obtain
\begin{align}
\nonumber
&\int_Q \left(|\chi F - \chi \nabla p_0 + J_1|^2 + |J_2|^2 
 + |\nabla_{x,t}(\nabla \chi \cdot u)|^2 \right) e^{2s\varphi} dxdt\\
\nonumber
\le &\; C\int_Q |F(x,t_0)|^2 e^{2s\varphi} dxdt 
+ C\int_Q |\nabla p_0|^2 e^{2s\varphi} dxdt \\
\label{stab:eq3}
& + C\int_Q \left(|J_1|^2 + |J_2|^2 
 + |\nabla_{x,t}(\nabla \chi \cdot u)|^2 \right) e^{2s\varphi} dxdt.
\end{align}
Moreover, by \eqref{stab:eq1}, we have
\begin{align*}
C\int_Q |\nabla p_0|^2 e^{2s\varphi} dxdt 
\le Ce^{Cs}\|p\|_{H^1(I;H^{\frac12}(\Gamma))}^2
\le Ce^{Cs}E_1^2.
\end{align*}
Here 
\begin{align*}
E_1 := \ 
&\|u\|_{H^2(I;L^2(\Gamma))} + \|\nabla u\|_{H^1(I;L^2(\Gamma))} 
+ \|H\|_{H^2(I;L^2(\Gamma))} + \|\nabla H\|_{H^1(I;L^2(\Gamma))} \\
&+ \|p\|_{H^1(I;H^{\frac12}(\Gamma))}.
\end{align*}
Further, from \eqref{con:chi}, we find the derivatives of $\chi$ vanish in 
$\{x\in \Omega; d(x)>\varepsilon_2\}$. 
Then we have
\begin{align*}
&C\int_Q \left(|J_1|^2 + |J_2|^2 
 + |\nabla_{x,t}(\nabla \chi \cdot u)|^2 \right) e^{2s\varphi} dxdt \\
\le &\; Ce^{2s\mu_1} \int_Q \left(|J_1|^2 + |J_2|^2 
 + |\nabla_{x,t}(\nabla \chi \cdot u)|^2 \right) dxdt\\
\le &\; Ce^{2s\mu_1} \left(\|u\|_{H^1(I;H^1(\Omega))}^2 
 + \|H\|_{L^2(I;H^1(\Omega))}^2 + \|p\|_{L^2(I;L^2(\Omega))}^2 \right) 
\le Ce^{2s\mu_1} M^2. 
\end{align*}
Here $M>0$ is the constant introduced in the statement of the theorem and  
\begin{equation}
\label{con:mu_1}
\mu_1 
:= \max\{\varphi; d(x)\le \varepsilon_2, t\in I\} 
= e^{\lambda \varepsilon_2}.
\end{equation}
Therefore, \eqref{stab:eq3} implies
\begin{align}
\nonumber
&\int_Q \left(|\chi F - \chi \nabla p_0 + J_1|^2 + |J_2|^2 
 + |\nabla_{x,t}(\nabla \chi \cdot u)|^2 \right) e^{2s\varphi} dxdt\\
\label{stab:eq4}
\le &\; C\int_Q |F(x,t_0)|^2 e^{2s\varphi} dxdt 
+ Ce^{Cs}E_1^2 + Ce^{2s\mu_1}M^2.
\end{align}

For the second integral on the right-hand side of \eqref{stab:eq2}, 
we note that $u_\chi, H_\chi$ and their derivatives vanish on 
$(\partial\Omega \setminus \Gamma) \times I$, thus we obtain
\begin{align}
\nonumber
& Cs^3\int_{\partial\Omega \times I} 
(|u_\chi|^2 + |H_\chi|^2 + |\nabla_{x,t}u_\chi|^2 + |\nabla_{x,t}H_\chi|^2) 
e^{2s\varphi} dSdt \\
\label{stab:eq5}
= &\; Cs^3\int_{\Gamma \times I} 
(|u_\chi|^2 + |H_\chi|^2 + |\nabla_{x,t}u_\chi|^2 + |\nabla_{x,t}H_\chi|^2) 
e^{2s\varphi} dSdt 
\le Cs^3e^{Cs} E_1^2.
\end{align}

For the third integral on the right-hand side of \eqref{stab:eq2}, we find that 
$\varphi(x,t) \le e^{\lambda (\|d\|_{C(\overline\Omega)} - \beta\delta^2)}$ 
for $x\in \Omega$ and $t=t_0 \pm \delta$. 
By denoting $\mu_2 := e^{\lambda (\|d\|_{C(\overline\Omega)} - \beta\delta^2)}$, 
we estimate
\begin{align}
\nonumber
&Cs^3\int_\Omega 
(|u_\chi|^2 + |H_\chi|^2 + |\nabla u_\chi|^2 + |\nabla H_\chi|^2) 
e^{2s\varphi} dx\Big|_{t=t_0\pm \delta} \\
\label{stab:eq6}
\le &\; Cs^3e^{2s\mu_2}\int_\Omega 
(|u|^2 + |H|^2 + |\nabla u|^2 + |\nabla H|^2) dx\Big|_{t=t_0\pm \delta} 
\le Cs^3e^{2s\mu_2} M^2.
\end{align}
Here we used the Sobolev embedding theorem in the last inequality. 
Therefore, inserting \eqref{stab:eq4}--\eqref{stab:eq6} into 
the right-hand side of \eqref{stab:eq2} yields
\begin{align}
\nonumber
&\int_Q \Biggl( 
\frac{1}{s^2}|\partial_t u_\chi|^2 + |\nabla u_\chi|^2 
+ s|\nabla H_\chi|^2 + s^2|u_\chi|^2 + s^3|H_\chi|^2
\Biggr) e^{2s\varphi} dxdt \\
\label{stab:eq7}
\le &\; C\int_Q |F(x,t_0)|^2 e^{2s\varphi} dxdt 
+ Cs^3e^{Cs}E_1^2 
+ Cs^3e^{2s\mu_3} M^2
\end{align}
for all sufficiently large $s\ge 1$. 
Here $\mu_3 := \max\{\mu_1,\mu_2\}$. 

Similarly, we employ Theorem \ref{thm:CE1} to \eqref{stab:sy2} and obtain
\begin{align}
\nonumber
&\int_Q \Biggl( 
\frac{1}{s^2}|\partial_t u_\chi^{(1)}|^2 + |\nabla u_\chi^{(1)}|^2 
+ s|\nabla H_\chi^{(1)}|^2 + s^2|u_\chi^{(1)}|^2 + s^3|H_\chi^{(1)}|^2
\Biggr) e^{2s\varphi} dxdt \\
\label{stab:eq8}
\le &\; C\int_Q |F(x,t_0)|^2 e^{2s\varphi} dxdt 
+ Cs^3e^{Cs}E_1^2 
+ Cs^3e^{2s\mu_3} M^2
\end{align}
for all sufficiently large $s\ge 1$.

\noindent {\bf Third Step: A preliminary estimate}

According to the first equation in \eqref{sy:MHD}, we see that
\begin{align*}
u^{(1)}_\chi(x,t_0) 
= &\; \chi\partial_t u(x,t_0)\\
= &\; \chi\left(\nu \Delta u(x,t_0) \!-\! (A^{(1)}\!\cdot\!\nabla)u(x,t_0)
 \!-\! (u\!\cdot\!\nabla)A^{(2)}(x,t_0) \!-\! \nabla p(x,t_0)\right)\\
& - \chi L_1 H(x,t_0) + \chi F(x,t_0), 
\end{align*}
and hence
\begin{align*}
\chi F(x,t_0) 
= & -\chi \Big( \nu \Delta u(x,t_0) - (A^{(1)}\!\cdot\!\nabla)u(x,t_0)
 - (u\!\cdot\!\nabla)A^{(2)}(x,t_0) - \nabla p(x,t_0)\Big) \\
& + u^{(1)}_\chi(x,t_0) + \chi L_1 H(x,t_0)
\end{align*}
for $x\in \Omega$. Recalling that $L_1 H$ contains at most first-order spatial derivatives of $H$ and 
$$
\chi L_1 H 
= L_1 H_\chi - (C^{(1)}\cdot \nabla\chi) H - \nabla\chi (C^{(3)}\cdot H),
$$
we estimate the weighted $L^2$-norm of $F(\cdot,t)$ in $\Omega$ at $t=t_0$: 
\begin{align}
\nonumber
\int_{\Omega} |\chi(x)F(x,t_0)|^2 e^{2s\varphi(x,t_0)} dx
\le &\; Ce^{Cs}E_2^2 
+ C\int_{\Omega} |u^{(1)}_\chi(x,t_0)|^2 e^{2s\varphi(x,t_0)} dx\\
\nonumber
& + C\int_{\Omega} (|\nabla H_\chi(x,t_0)|^2 
 + |H_\chi(x,t_0)|^2) e^{2s\varphi(x,t_0)} dx\\
\label{stab:eq9}
& + C\int_{\Omega} |\nabla\chi(x)|^2 |H(x,t_0)|^2 e^{2s\varphi(x,t_0)} dx.
\end{align}
Here 
$$
E_2 := 
\|u(\cdot,t_0)\|_{H^2(\Omega)} 
+ \|\nabla p(\cdot,t_0)\|_{L^2(\Omega)}.
$$
Next we estimate the second, the third and the fourth terms on the right-hand side above respectively. 
For the second term, we have
\begin{align*}
&\int_{\Omega} |u_\chi^{(1)}(x,t_0)|^2 e^{2s\varphi(x,t_0)} dx\\
=&\; \int^{t_0}_{t_0-\delta} 
 \left( \partial_t \int_{\Omega} |u_\chi^{(1)}|^2 e^{2s\varphi} dx\right) dt 
 + \int_{\Omega} |u_\chi^{(1)}(x,t_0-\delta)|^2 e^{2s\varphi(x,t_0-\delta)} dx\\
=&\; \int^{t_0}_{t_0-\delta} \int_{\Omega} \left(2u_\chi^{(1)}\!\cdot\!\partial_t u_\chi^{(1)} 
 + 2s(\partial_t \varphi)| u_\chi^{(1)}|^2 \right) e^{2s\varphi} dxdt 
 + \int_{\Omega} |u_\chi^{(1)}(x,t_0-\delta)|^2 e^{2s\varphi(x,t_0-\delta)} dx\\
\le &\; C\int_{Q} \left(|u_\chi^{(1)}| |\partial_t u_\chi^{(1)}| 
 + s|u_\chi^{(1)}|^2 \right) e^{2s\varphi} dxdt 
 + \int_{\Omega} |u_\chi^{(1)}(x,t_0-\delta)|^2 e^{2s\varphi(x,t_0-\delta)} dx\\
\le &\; C\int_Q \left( \frac{1}{s^2} |\partial_t u_\chi^{(1)}|^2 
 + s^2|u_\chi^{(1)}|^2 \right) e^{2s\varphi} dxdt 
 + Ce^{2s\mu_3} M^2.
\end{align*}
Here we used $|a||b| \le \frac{1}{2}(s^2|a|^2 + s^{-2}|b|^2)$ and 
the argument used in \eqref{stab:eq6}. 
Similarly, we have
\begin{align*}
&\int_{\Omega} |H_\chi(x,t_0)|^2 e^{2s\varphi(x,t_0)} dx \\
&\le C\int_{Q} \left( \frac{1}{s^2} |\partial_t H_\chi|^2 
 + s^2|H_\chi|^2 \right) e^{2s\varphi} dxdt 
 + \int_{\Omega} |H_\chi(x,t_0-\delta)|^2 e^{2s\varphi(x,t_0-\delta)} dx\\
&\le C\int_{Q} \left( \frac{1}{s^2} |H_\chi^{(1)}|^2 
 + s^2|H_\chi|^2 \right) e^{2s\varphi} dxdt
 + Ce^{2s\mu_3} M^2.
\end{align*}
Moreover, we calculate for $j=1,2,3$,
\begin{align*}
&\int_{\Omega} |\partial_j H_\chi(x,t_0)|^2 e^{2s\varphi(x,t_0)} dx\\
=&\; \int^{t_0}_{t_0-\delta} 
 \left( \partial_t \int_{\Omega} |\partial_j H_\chi|^2 e^{2s\varphi} dx\right) dt 
 + \int_{\Omega} |\partial_j H_\chi(x,t_0-\delta)|^2 e^{2s\varphi(x,t_0-\delta)} dx\\
\le &\; C\int_{Q} (|\partial_j H_\chi| |\partial_j \partial_t H_\chi| 
 + s|\partial_j H_\chi|^2) e^{2s\varphi} dxdt 
 + Ce^{2s\mu_3} M^2\\
\le &\; C\int_{Q} \left( \frac{1}{s} |\partial_j H_\chi^{(1)}|^2 
 + s|\partial_j H_\chi|^2 \right) e^{2s\varphi} dxdt 
 + Ce^{2s\mu_3} M^2.
\end{align*}
For the fourth term, we note again that the derivatives of $\chi$ vanish in 
$\{x\in \Omega; d(x)>\varepsilon_2\}$, which implies 
\begin{align*}
&\int_{\Omega} |\nabla\chi(x)|^2 |H(x,t_0)|^2 e^{2s\varphi(x,t_0)} dx 
\le Ce^{2s\mu_1} \|H(\cdot,t_0)\|_{L^2(\Omega)}^2
\le Ce^{2s\mu_3} M^2.
\end{align*}
Therefore, equipped with all the above estimations, \eqref{stab:eq9} together with 
\eqref{stab:eq7}--\eqref{stab:eq8} implies
\begin{equation}
\label{stab:eq10}
\int_{\Omega} |\chi(x)F(x,t_0)|^2 e^{2s\varphi(x,t_0)} dx
\le C\int_Q |F(x,t_0)|^2 e^{2s\varphi} dxdt 
 + Cs^3e^{2s\mu_3} M^2 
 + Cs^3e^{Cs}D_1^2
\end{equation}
for all sufficiently large $s\ge 1$.  
Here we noted that $E_1^2 + E_2^2 \le C D_1^2$ where 
$D_1$ is the norm of the data introduced in the theorem. 

\noindent {\bf Fourth step: Absorption and completion of the proof}

Now we absorb the first term on the right-hand side of \eqref{stab:eq10} into the left-hand side. 
With reference to $\varphi(x,t) = d(x) - \beta (t-t_0)^2$, $(x,t)\in Q$, 
the Lebesgue's dominated convergence theorem yields 
\begin{align*}
\int_Q |F(x,t_0)|^2 e^{2s\varphi} dxdt
=& \; o(1) \int_{\Omega} |F(x,t_0)|^2 e^{2s\varphi(x,t_0)} dx
\end{align*}
as $s \to \infty$. 
We refer to \cite[Lemma 5.3, pp.120]{BY17} for the detailed argument. 
By recalling the choice of $\mu_1$, see \eqref{con:mu_1} and that $F\in \mathcal F_{M,1}$, we have
\begin{align*}
\quad& \int_Q |F(x,t_0)|^2 e^{2s\varphi} dxdt \\
=&\; o(1) \left( \int_{\{x\in\Omega;\varphi(x,t_0)>\mu_1\}}
 + \int_{\{x\in\Omega;\varphi(x,t_0)\le\mu_1\}} \right) |F(x,t_0)|^2 e^{2s\varphi(x,t_0)} dx\\
=&\; o(1) \int_{\{x\in\Omega;\varphi(x,t_0)>\mu_1\}} |F(x,t_0)|^2 e^{2s\varphi(x,t_0)} dx
 + Ce^{2s\mu_1}M^2.
\end{align*}
On the other hand, with reference to \eqref{con:chi}, we have $\chi = 1$ in 
$\{x\in \Omega; d(x)>\varepsilon_2\} = \{x\in \Omega; \varphi(x,t_0) > \mu_1\}$. 
Then
\begin{align*}
\int_{\Omega} |\chi(x)F(x,t_0)|^2 e^{2s\varphi(x,t_0)} dx
\ge &\; \int_{\{x\in\Omega;\varphi(x,t_0)>\mu_1\}} |\chi(x)|^2 |F(x,t_0)|^2 e^{2s\varphi(x,t_0)} dx\\
= &\; \int_{\{x\in\Omega;\varphi(x,t_0)>\mu_1\}} |F(x,t_0)|^2 e^{2s\varphi(x,t_0)} dx,
\end{align*}
so that \eqref{stab:eq10} leads to
\begin{align*}
&\int_{\{x\in\Omega;\varphi(x,t_0)>\mu_1\}} |F(x,t_0)|^2 e^{2s\varphi(x,t_0)} dx\\
\le &\; o(1) \int_{\{x\in\Omega;\varphi(x,t_0)>\mu_1\}} |F(x,t_0)|^2 e^{2s\varphi(x,t_0)} dx
 + Cs^3e^{2s\mu_3}M^2 
 + Cs^3e^{Cs}D_1^2.
\end{align*}
Absorbing the first term on the right-hand side into the left-hand side, we obtain
$$
\int_{\{x\in\Omega;\varphi(x,t_0)>\mu_1\}} |F(x,t_0)|^2 e^{2s\varphi(x,t_0)} dx
\le Cs^3e^{2s\mu_3}M^2 + Cs^3e^{Cs}D_1^2
$$
for sufficiently large $s$. 
By noting that $\mu_1=e^{\lambda \varepsilon_2}$, $0<\varepsilon_2<\varepsilon_3$ 
and denoting $\mu_0 := e^{\lambda \varepsilon_3}$, 
we give a lower bound of the left-hand side above: 
\begin{align*}
\int_{\{x\in\Omega;\varphi(x,t_0)>\mu_1\}} |F(x,t_0)|^2 e^{2s\varphi(x,t_0)} dx 
&= \int_{\{x\in\Omega;d(x)>\varepsilon_2\}} |F(x,t_0)|^2 e^{2s\varphi(x,t_0)} dx\\
&\ge \int_{\{x\in\Omega;d(x)>\varepsilon_3\}} |F(x,t_0)|^2 e^{2s\varphi(x,t_0)} dx\\
&> e^{2s\mu_0}\int_{\{x\in\Omega;d(x)>\varepsilon_3\}} |F(x,t_0)|^2 dx\\
&\ge e^{2s\mu_0}\|F(\cdot,t_0)\|_{L^2(\Omega_0)}^2
\end{align*}
where we used the inclusion 
$\Omega_0 \subset \{x\in \Omega; d(x)>\varepsilon_3\}$ 
in the last inequality. 
Therefore, 
\begin{equation}
\label{stab:eq11}
\|F(\cdot,t_0)\|^2_{L^2(\Omega_0)} 
\le Cs^3e^{-2s(\mu_0-\mu_3)}M^2 + Cs^3e^{Cs}D_1^2
\end{equation}
for sufficiently large $s$. 
Finally, by taking a large $\beta>0$ such that
$$
\varepsilon_3 > \|d\|_{C(\overline\Omega)} - \beta\delta^2, 
$$
which implies $\mu_0>\mu_2$, and noting $\mu_0>\mu_1$ from the definition, 
we find $\mu_0>\max\{\mu_1,\mu_2\}=\mu_3$. 
Then we end up with the stability of H\"older type by minimizing 
the right-hand side of \eqref{stab:eq11} with respect to $s$ 
(see e.g., \cite[Proof of Theorem 5.1, pp. 27--28]{Y09}).
\subsection{Proof of Theorem \ref{thm:stab2}}
\label{subsec:proof2}

With the additional regularity assumptions on the unknown and the measurements, 
we can employ Theorem \ref{thm:CE2} and derive the stability result without 
the measurements on the pressure $p$. 

\noindent {\bf First step: A preliminary estimate}

We start the proof with a weighted estimate. 
In order to erase the effect of the pressure $p$, 
we need to apply the rotation operator 
to both sides of the first equation of \eqref{sy:MHD}. 
Actually by using the equalities 
$\mathrm{rot} (\mathrm{rot}\, u) = - \Delta u + \nabla (\mathrm{div}\, u)$, 
$\mathrm{rot} (\mathrm{rot}\, F) = - \Delta F + \nabla (\mathrm{div}\, F)$
and $\mathrm{rot}(\nabla p) = 0$, 
we take the rotation operator twice on both sides of the first equation of \eqref{sy:MHD}. 
In particular, we fix $t=t_0$ and then by \eqref{con:F2-2b} we obtain
\begin{align}
\label{stab2:eq1}
\partial_t \Delta u(x,t_0) - \nu\Delta^2 u(x,t_0) - \mathrm{rot}(\mathrm{rot}\, L_1 H(x,t_0)) 
= \Delta F(x,t_0) + L_3 u(x,t_0),                                                                    \quad x\in \Omega
\end{align}
where $L_3$ is a third-order differential operator given by 
$$
L_3 u := \mathrm{rot}\left(\mathrm{rot}\left((A^{(1)}\cdot\nabla) u + (u\cdot\nabla) A^{(2)}\right)\right). 
$$
By taking the weighted $L^2$-norm over $\Omega$ on both sides of \eqref{stab2:eq1}, we obtain
\begin{align}
\nonumber
&\int_\Omega |\Delta F(x,t_0)|^2 e^{2s\varphi(x,t_0)}dx \\
\nonumber
= &\; \int_\Omega |\partial_t \Delta u(x,t_0) - \nu\Delta^2 u(x,t_0) 
 - \mathrm{rot}(\mathrm{rot}\, L_1 H(x,t_0)) - L_3 u(x,t_0)|^2 e^{2s\varphi(x,t_0)} dx \\
\label{stab2:eq2}
\le &\; C\int_\Omega |\partial_t \Delta u(x,t_0)|^2 e^{2s\varphi(x,t_0)} dx 
 + Ce^{Cs}E_3^2
\end{align}
where 
\begin{align*}
E_3 := \|u(\cdot,t_0)\|_{H^4(\Omega)} + \|H(\cdot,t_0)\|_{H^3(\Omega)}.
\end{align*}
Next, we give an upper bound of the first term on the right-hand side of \eqref{stab2:eq2}. 

Actually, similarly as we have done in the former subsection, we estimate
\begin{align}
\nonumber
&\int_{\Omega} |\partial_t \Delta u(x,t_0)|^2 e^{2s\varphi(x,t_0)}dx\\
\nonumber
= &\; \int_{t_0-\delta}^{t_0} \partial_t 
 \left(\int_\Omega |\partial_t \Delta u|^2 e^{2s\varphi}dx \right) dt
 + \int_\Omega |\partial_t \Delta u(x,t_0-\delta)|^2 e^{2s\varphi(x,t_0-\delta)} dx \\
\nonumber
= &\; \int_{t_0-\delta}^{t_0} \int_\Omega \left( 2\Delta (\partial_t u)\cdot \Delta (\partial_t^2 u) 
 + 2s(\partial_t \varphi)|\Delta (\partial_t u)|^2 \right)e^{2s\varphi} dxdt\\
\nonumber
& + \int_\Omega |\partial_t \Delta u(x,t_0-\delta)|^2 e^{2s\varphi(x,t_0-\delta)} dx \\
\nonumber
\le &\; C\int_Q \left( |\Delta (\partial_t u)| |\Delta (\partial_t^2 u)| 
 + s|\Delta (\partial_t u)|^2 \right)e^{2s\varphi} dxdt\\
\nonumber
& + \int_\Omega |\partial_t \Delta u(x,t_0-\delta)|^2 e^{2s\varphi(x,t_0-\delta)} dx \\
\label{stab2:eq3}
\le &\; C\int_Q \left(|\Delta (\partial_t^2 u)|^2 
 + s|\Delta (\partial_t u)|^2 \right)e^{2s\varphi} dxdt 
 + Ce^{2s\mu_1}M^2
\end{align}
where $M>0$ is the constant introduced in the statement of the theorem 
and $\mu_1$ is the maximum of $\varphi$ over $\overline\Omega \times \{ t_0 \pm \delta\}$ 
and $(\partial\Omega \setminus \Gamma) \times I$, which implies
\begin{align}
\nonumber
&\mu_1 \ge \max\{\varphi(x,t); x\in \overline\Omega, t = t_0 \pm \delta\} 
= e^{\lambda (\|d\|_{C(\overline\Omega)} - \beta \delta^2)} \quad \mbox{ and }\\
\label{con2:mu_1}
&\mu_1 \ge \max\{\varphi(x,t); x\in \partial\Omega \setminus \Gamma, t \in I\} = 1.
\end{align}
Here we used $2A\cdot B \le |A|^2 + |B|^2$ and \eqref{stab2:eq3} holds true for all $s\ge 1$. 
Thus, we reach
\begin{align}
\nonumber
\int_\Omega |\Delta F(x,t_0)|^2 e^{2s\varphi(x,t_0)}dx 
\le &\; C\int_Q \left(|\Delta (\partial_t^2 u)|^2 
 + s|\Delta (\partial_t u)|^2 \right)e^{2s\varphi} dxdt \\
\label{stab2:eq4}
& + Ce^{2s\mu_1}M^2 + Ce^{Cs}E_3^2. 
\end{align}

\noindent{\bf Second step: Application of the Carleman estimate} 

For simplicity of the expressions, we introduce the following notations: 
$$
u^{(k)} = \partial_t^k u, \ 
H^{(k)} = \partial_t^k H, \ 
p^{(k)} = \partial_t^k p,  \quad k=1,2.
$$
Then taking once and twice time derivatives of the governing system \eqref{sy:MHD} yields 
\begin{equation}
\label{sy:MHD1}
\left\{
\begin{aligned}
&\partial_t u^{(1)} \!-\! \nu\Delta u^{(1)} + (A^{(1)}\!\cdot\!\nabla) u^{(1)} 
 + (u^{(1)}\!\cdot\!\nabla) A^{(2)} + L_1 H^{(1)} + \nabla p^{(1)} = \partial_t F + J_1, \\
&\partial_t H^{(1)} - \kappa\Delta H^{(1)} + (A^{(3)}\cdot\nabla) H^{(1)} 
 + (H^{(1)}\cdot\nabla) A^{(4)} + L_2 u^{(1)} = J_2, \\
& \mathrm{div}\, u^{(1)} = 0,
\end{aligned}
\right.
\end{equation}
and
\begin{equation}
\label{sy:MHD2}
\left\{
\begin{aligned}
&\partial_t u^{(2)} \!-\! \nu\Delta u^{(2)} + (A^{(1)}\!\cdot\!\nabla) u^{(2)} 
 + (u^{(2)}\!\cdot\!\nabla) A^{(2)} \!+\! L_1 H^{(2)} + \nabla p^{(2)} = \partial_t^2 F + J_3,\\
&\partial_t H^{(2)} - \kappa\Delta H^{(2)} + (A^{(3)}\cdot\nabla) H^{(2)} 
 + (H^{(2)}\cdot\nabla) A^{(4)} + L_2 u^{(2)} = J_4, \\
& \mathrm{div}\, u^{(2)} = 0 
\end{aligned}
\right.
\end{equation}
Here $J_1, J_2$ denote the terms including at most first-order spatial derivative of $u$ and $H$ 
while $J_3, J_4$ include additionally at most first-order spatial derivatives of 
$u^{(1)}, H^{(1)}$, 
and we have 
\begin{align*}
|\mathrm{rot}\, J_k| 
\le &\; C\left(|u| + |\nabla u| + |\nabla \mathrm{rot}\, u| 
 + |H| + |\nabla H| + |\nabla \mathrm{rot}\, H|\right),                              \quad k=1,2. \\
|\mathrm{rot}\, J_k| 
\le &\; C\Big(|u| + |\nabla u| + |\nabla \mathrm{rot}\, u| 
 + |H| + |\nabla H| + |\nabla \mathrm{rot}\, H| + |u^{(1)}| + |\nabla u^{(1)}| \\
&\quad + |\nabla \mathrm{rot}\, u^{(1)}| + |H^{(1)}| 
 + |\nabla H^{(1)}| + |\nabla \mathrm{rot}\, H^{(1)}|\Big),                       \quad k=3,4.
\end{align*}
We employ Theorem \ref{thm:CE2} to \eqref{sy:MHD}, and in particular, we obtain 
\begin{align}
\nonumber
&\int_Q \left(s|\nabla \mathrm{rot}\, u|^2 + s|\nabla \mathrm{rot}\, H|^2 
 + s^2|\nabla u|^2 + s^2|\nabla H|^2 + s^4|u|^2 + s^4|H|^2 \right)e^{2s\varphi} dxdt\\
\nonumber
&\le C\int_Q |\mathrm{rot}\, F|^2 e^{2s\varphi} dxdt \\
\nonumber
& + Cs^4\int_{\partial\Omega \times I} \left(|\nabla_{x,t} \mathrm{rot}\, u|^2 
 + |\nabla u|^2 + |\nabla_{x,t} \mathrm{rot}\, H|^2 
 + |\nabla H|^2 + |u|^2 + |H|^2 \right)e^{2s\varphi} dSdt\\
\label{stab2:eq5}
& + Cs^3\int_\Omega \left(|\nabla \mathrm{rot}\, u|^2 
 + |\mathrm{rot}\, u|^2 + |\nabla \mathrm{rot}\, H|^2 
 + |\mathrm{rot}\, H|^2 \right)e^{2s\varphi} dx \Big|_{t=t_0\pm \delta} 
\end{align}
for all sufficiently large $s\ge 1$. 
Now we estimate the upper bounds of the second and the third terms 
on the right-hand side above. 

By noting $\varphi(\cdot,t_0\pm \delta) \le \mu_1$ on $\overline\Omega$ 
($\mu_1$ is given in \eqref{con2:mu_1}), the third term admits 
\begin{align*}
Cs^3\int_\Omega \left(|\nabla \mathrm{rot}\, u|^2 
 + |\mathrm{rot}\, u|^2 + |\nabla \mathrm{rot}\, H|^2 
 + |\mathrm{rot}\, H|^2 \right)e^{2s\varphi} dx \Big|_{t=t_0\pm \delta} 
\le Cs^3e^{2s\mu_1}M^2. 
\end{align*}
For the second term, we divide the integral into two parts over $\Gamma \times I$ and 
$(\partial\Omega\setminus \Gamma) \times I$. 
By noting that $\varphi$ is dominated by some constant in $\Gamma \times I$, 
and is dominated by $\mu_1\ge 1$ in $(\partial\Omega\setminus \Gamma) \times I$, 
we obtain
\begin{align*}
&Cs^4\int_{\partial\Omega \times I} \left(|\nabla_{x,t} \mathrm{rot}\, u|^2 
 + |\nabla u|^2 + |\nabla_{x,t} \mathrm{rot}\, H|^2 + |\nabla H|^2 
 + |u|^2 + |H|^2 \right)e^{2s\varphi} dSdt \\
&= Cs^4\int_{\Gamma \times I} \left(|\nabla_{x,t} \mathrm{rot}\, u|^2 
 + |\nabla u|^2 + |\nabla_{x,t} \mathrm{rot}\, H|^2 + |\nabla H|^2 
 + |u|^2 + |H|^2 \right)e^{2s\varphi} dSdt \\
&+ Cs^4\int_{(\partial\Omega\setminus \Gamma) \times I} 
 \left(|\nabla_{x,t} \mathrm{rot}\, u|^2 + |\nabla u|^2 + |\nabla_{x,t} \mathrm{rot}\, H|^2 
 + |\nabla H|^2 + |u|^2 + |H|^2 \right)e^{2s\varphi} dSdt \\
&\le Cs^4 e^{Cs}E_4^2 + Cs^4e^{2s\mu_1}M^2
\end{align*}
where
\begin{align*}
E_4 := &\; 
\|\nabla \mathrm{rot}\, u\|_{H^2(I;L^2(\Gamma))} 
+ \|\nabla \mathrm{rot}\, H\|_{H^2(I;L^2(\Gamma))} 
+ \|\nabla u\|_{H^3(I;L^2(\Gamma))} \\
&+ \|\nabla H\|_{H^3(I;L^2(\Gamma))} 
+ \|u\|_{H^2(I;L^2(\Gamma))} 
+ \|H\|_{H^2(I;L^2(\Gamma))}.
\end{align*}
Therefore, we derive from \eqref{stab2:eq5} that
\begin{align}
\nonumber
&\int_Q \left(s|\nabla \mathrm{rot}\, u|^2 + s|\nabla \mathrm{rot}\, H|^2 
 + s^2|\nabla u|^2 + s^2|\nabla H|^2 + s^4|u|^2 + s^4|H|^2 \right)e^{2s\varphi} dxdt\\
\label{stab2:eq6}
&\le C\int_Q |\mathrm{rot}\, F|^2 e^{2s\varphi} dxdt 
 + Cs^4 e^{Cs}E_4^2 + Cs^4e^{2s\mu_1}M^2
\end{align}
for all sufficiently large $s\ge 1$. 
Similarly, we employ Theorem \ref{thm:CE2} to \eqref{sy:MHD1} and we obtain
\begin{align}
\nonumber
&\int_Q \Big(s|\Delta u^{(1)}|^2 + s|\nabla \mathrm{rot}\, u^{(1)}|^2 
 + s|\nabla \mathrm{rot}\, H^{(1)}|^2 + s^2|\nabla u^{(1)}|^2  \\
\nonumber
&\qquad + s^2|\nabla H^{(1)}|^2 + s^4|u^{(1)}|^2 + s^4|H^{(1)}|^2 \Big)e^{2s\varphi} dxdt\\
\nonumber
&\le C\int_Q \left(|\mathrm{rot}(\partial_t F) + \mathrm{rot}\, J_1|^2 
 + |\mathrm{rot}\, J_2|^2\right) e^{2s\varphi} dxdt 
 + Cs^4 e^{Cs}E_4^2 + Cs^4e^{2s\mu_1}M^2\\
\nonumber
&\le C\int_Q |\mathrm{rot}(\partial_t F)|^2 e^{2s\varphi} dxdt 
 + Cs^4 e^{Cs}E_4^2 + Cs^4e^{2s\mu_1}M^2\\
\label{stab2:eq7}
&+ C \int_Q \left(|\nabla \mathrm{rot}\, u|^2 + |\nabla \mathrm{rot}\, H|^2 
 + |\nabla u|^2 + |\nabla H|^2 + |u|^2 + |H|^2 \right)e^{2s\varphi} dxdt
\end{align}
for all sufficiently large $s\ge 1$. 
Inserting \eqref{stab2:eq6} into the right-hand side of \eqref{stab2:eq7} yields
\begin{align}
\nonumber
&\int_Q \Big(s|\Delta u^{(1)}|^2 + s|\nabla \mathrm{rot}\, u^{(1)}|^2 
 + s|\nabla \mathrm{rot}\, H^{(1)}|^2 + s^2|\nabla u^{(1)}|^2  \\
\nonumber
&\qquad + s^2|\nabla H^{(1)}|^2 + s^4|u^{(1)}|^2 
 + s^4|H^{(1)}|^2 \Big)e^{2s\varphi} dxdt\\
\label{stab2:eq8}
&\le C\int_Q \left(|\mathrm{rot}(\partial_t F)|^2 
 + |\mathrm{rot}\, F|^2\right)e^{2s\varphi} dxdt 
 + Cs^4 e^{Cs}E_4^2 + Cs^4e^{2s\mu_1}M^2
\end{align}
for all sufficiently large $s\ge 1$. 

Next we employ Theorem \ref{thm:CE2} to \eqref{sy:MHD2}. 
We pick up the necessary terms and obtain
\begin{align}
\nonumber
&\int_Q s|\Delta u^{(2)}|^2 e^{2s\varphi} dxdt\\
\nonumber
&\le C\int_Q \left(|\mathrm{rot}(\partial_t^2 F) + \mathrm{rot}\, J_3|^2 
 + |\mathrm{rot}\, J_4|^2\right) e^{2s\varphi} dxdt 
 + Cs^4 e^{Cs}E_4^2 + Cs^4e^{2s\mu_1}M^2\\
\nonumber
&\le C\int_Q |\mathrm{rot}(\partial_t^2 F)|^2 e^{2s\varphi} dxdt 
 + Cs^4 e^{Cs}E_4^2 + Cs^4e^{2s\mu_1}M^2\\
\nonumber
&+ C \int_Q \Big(|\nabla \mathrm{rot}\, u|^2 + |\nabla \mathrm{rot}\, H|^2 
 + |\nabla u|^2 + |\nabla H|^2 + |u|^2 + |H|^2 + |\nabla \mathrm{rot}\, u^{(1)}|^2 \\
\label{stab2:eq9}
&+ |\nabla \mathrm{rot}\, H^{(1)}|^2 + |\nabla u^{(1)}|^2 + |\nabla H^{(1)}|^2 
 + |u^{(1)}|^2 + |H^{(1)}|^2 \Big)e^{2s\varphi} dxdt
\end{align}
for all sufficiently large $s\ge 1$. 
By \eqref{stab2:eq6} and \eqref{stab2:eq8}, 
we can get rid of the last term on the right-hand side of \eqref{stab2:eq9}. 
Moreover, by \eqref{con:F2-2a} we have
$$
|\mathrm{rot}(\partial_t^k F)(x,t)|^2 
\le C\left(|\nabla F(x,t_0)|^2 + |F(x,t_0)|^2\right),                      \quad k=0,1,2,
$$
then we reach
\begin{align*}
&\int_Q \left(s|\Delta u^{(1)}|^2 + s|\Delta u^{(2)}|^2 \right)e^{2s\varphi} dxdt\\
&\le C\int_Q \left(|\nabla F(x,t_0)|^2 + |F(x,t_0)|^2\right)e^{2s\varphi} dxdt 
 + Cs^4 e^{Cs}E_4^2 + Cs^4e^{2s\mu_1}M^2
\end{align*}
for all sufficiently large $s\ge 1$. 
We insert this inequality into the right-hand side of \eqref{stab2:eq4} and we obtain
\begin{align}
\nonumber
\int_\Omega |\Delta F(x,t_0)|^2 e^{2s\varphi(x,t_0)}dx 
\le &\; C\int_Q \left(|\nabla F(x,t_0)|^2 + |F(x,t_0)|^2\right)e^{2s\varphi} dxdt \\
\label{stab2:eq9.5}
& + Cs^4e^{2s\mu_1}M^2 + Cs^4e^{Cs}D_2^2
\end{align}
for all sufficiently large $s\ge 1$. 
Here we noted that $E_3^2 + E_4^2 \le C D_2^2$ where 
$D_2$ is the norm of the data that introduced in the theorem. 

\noindent {\bf Third Step: Lower bound of \eqref{stab2:eq9.5} and completion of the proof}

We introduce the following lemma, which is based on a classical Carleman estimate for 
elliptic equations. 
\begin{lem}
\label{lem:CE-ellip}
There exist constants $\widehat s>0$ and $C>0$ such that 
\begin{align}
\nonumber
&\quad \int_\Omega \left(s|\nabla F(x,t_0)|^2 
 + s^3|F(x,t_0)|^2\right)e^{2s\varphi(x,t_0)} dx \\
\label{stab2:eq10}
&\le C\int_\Omega |\Delta F(x,t_0)|^2 e^{2s\varphi(x,t_0)}dx 
 + Cs^3e^{2s\mu_1}M^2
\end{align}
for all $s\ge \widehat s$. 
\end{lem}
\begin{proof}
According to the assumption that $F\in \mathcal F_{2,M}$, 
we have $F(\cdot,t_0)$ and $\nabla F(\cdot,t_0)$ vanish on $\Gamma$. 
By the classical Carleman estimate for the elliptic equations, we obtain 
\begin{align*}
&\int_\Omega \left(s|\nabla F(x,t_0)|^2 
 + s^3|F(x,t_0)|^2\right)e^{2s\varphi(x,t_0)} dx \\
&\le C\int_\Omega |\Delta F(x,t_0)|^2 e^{2s\varphi(x,t_0)}dx 
+ C\int_{\partial\Omega\setminus \Gamma} (s|\nabla F(x,t_0)|^2 
+ s^3|F(x,t_0)|^2) e^{2s\varphi(x,t_0)}dx
\end{align*}
for all sufficiently large $s\ge 1$. 
By noting 
$\varphi(\cdot,t_0) = e^{\lambda d} = 1 \le \mu_1$ 
on $\partial\Omega \setminus \Gamma$, 
we have
\begin{align*}
&\int_\Omega \left(s|\nabla F(x,t_0)|^2 
 + s^3|F(x,t_0)|^2\right)e^{2s\varphi(x,t_0)} dx \\
&\le C\int_\Omega |\Delta F(x,t_0)|^2 e^{2s\varphi(x,t_0)}dx 
 + Cs^3e^{2s\mu_1}M^2
\end{align*}
for all sufficiently large $s\ge 1$. 
\end{proof}

Now we combine \eqref{stab2:eq9.5} with \eqref{stab2:eq10} and we obtain
\begin{align}
\nonumber
&\int_\Omega \left(s|\nabla F(x,t_0)|^2 
 + s^3|F(x,t_0)|^2\right)e^{2s\varphi(x,t_0)} dx \\
\label{stab2:eq11}
&\le C\int_Q \left(|\nabla F(x,t_0)|^2 
 + |F(x,t_0)|^2\right)e^{2s\varphi} dxdt 
 + Cs^4e^{2s\mu_1}M^2 + Cs^4e^{Cs}D_2^2
\end{align}
for all sufficiently large $s\ge 1$. 
Since $\varphi(x,\cdot)$ attains its maximum at $t_0$, we have
$$
\int_Q \left(|\nabla F(x,t_0)|^2 + |F(x,t_0)|^2\right)e^{2s\varphi} dxdt 
\le 2\delta \int_\Omega \left(|\nabla F(x,t_0)|^2 + |F(x,t_0)|^2\right)e^{2s\varphi(x,t_0)} dx.
$$
Thus, by taking $s$ large enough, we can absorb the first term on 
the right-hand side of \eqref{stab2:eq11} into the left-hand side. 
This yields
\begin{align}
\label{stab2:eq12}
\int_\Omega \left(s|\nabla F(x,t_0)|^2 + s^3|F(x,t_0)|^2\right)e^{2s\varphi(x,t_0)} dx 
\le &\; Cs^4e^{2s\mu_1}M^2 + Cs^4e^{Cs}D_2^2
\end{align}
for all $s$ large enough. 

We recall that $\overline{\Omega_0} \subset \Omega \cup \Gamma$, 
which implies that there exists $\varepsilon_0>0$ such that 
$d(x)\ge \varepsilon_0$ for $x\in \Omega_0$, that is, 
$\varphi(x,t_0) = e^{\lambda d(x)} \ge e^{\lambda \varepsilon_0} =: \mu_0$. 
Therefore, \eqref{stab2:eq12} yields
\begin{align*}
e^{2s\mu_0}\|F(\cdot,t_0)\|_{H^1(\Omega_0)}^2
&\le \int_{\Omega_0} \left(s|\nabla F(x,t_0)|^2 
 + s^3|F(x,t_0)|^2\right)e^{2s\varphi(x,t_0)} dx \\
&\le \int_{\Omega} \left(s|\nabla F(x,t_0)|^2 
 + s^3|F(x,t_0)|^2\right)e^{2s\varphi(x,t_0)} dx \\
&\le Cs^4e^{2s\mu_1}M^2 + Cs^4e^{Cs}D_2^2, 
\end{align*}
and hence
\begin{align}
\label{stab2:eq13}
\|F(\cdot,t_0)\|_{H^1(\Omega_0)}^2 
\le Cs^4e^{-2s(\mu_0-\mu_1)}M^2 + Cs^4e^{Cs}D_2^2
\end{align}
for all $s$ large enough. 
Since $\mu_1 = \max\{e^{\lambda(\|d\|_{C(\overline\Omega)} - \beta\delta^2)}, 1\}$ 
and $\mu_0 = e^{\lambda \varepsilon_0} >1$, 
we can find $\beta>0$ large enough such that 
$\|d\|_{C(\overline\Omega)} - \beta\delta^2 < \varepsilon_0$ 
and consequently $\mu_0-\mu_1>0$. 
This enables us to end up with a stability of H\"older type by minimizing 
the right-hand side of \eqref{stab2:eq13} with respect to $s$. 
\subsection{Proof of Theorem \ref{thm:stab3}}

Thanks to the assumption \eqref{con:F2-1}, we need only apply the rotation operator once to 
the governing equation \eqref{sy:MHD}. 
Actually we can simply follow the steps of the former subsection. 

\noindent {\bf First step: A preliminary estimate}

Again we start the proof with a weighted estimate. 
In order to erase the effect of the pressure $p$, 
we apply the rotation operator to both sides of the first equation of \eqref{sy:MHD}. 
In particular, by setting $v=\mathrm{rot}\, u$ and taking $t=t_0$, we obtain
\begin{align}
\label{stab3:eq1}
\partial_t v(x,t_0) - \nu \Delta v(x,t_0) + \mathrm{rot}\, L_1 H(x,t_0)
= \mathrm{rot}\, F(x,t_0) + L_4 u(x,t_0),                           \quad x\in \Omega,
\end{align}
where $L_4$ is a second-order differential operator given by 
$$
L_4 u := -\mathrm{rot}\left((A^{(1)}\cdot\nabla) u + (u\cdot\nabla) A^{(2)}\right). 
$$
By taking the weighted $L^2$-norm over $\Omega$ on both sides of \eqref{stab3:eq1}, we obtain
\begin{align}
\nonumber
&\int_\Omega |\mathrm{rot}\, F(x,t_0)|^2 e^{2s\varphi(x,t_0)}dx \\
\nonumber
= &\; \int_\Omega |\partial_t v(x,t_0) - \nu\Delta v(x,t_0) 
 + \mathrm{rot}\, L_1 H(x,t_0) - L_4 u(x,t_0)|^2 e^{2s\varphi(x,t_0)} dx \\
\nonumber
\le &\; C\int_\Omega |\partial_t v(x,t_0)|^2 e^{2s\varphi(x,t_0)} dx 
 + C\int_\Omega |\nabla \mathrm{rot}\, H(x,t_0)|^2 e^{2s\varphi(x,t_0)} dx\\
\label{stab3:eq2}
&+ C\int_\Omega |\nabla H(x,t_0)|^2 e^{2s\varphi(x,t_0)} dx 
 + C\int_\Omega |H(x,t_0)|^2 e^{2s\varphi(x,t_0)} dx 
 + Ce^{Cs}E_5^2.
\end{align}
Here
\begin{align*}
E_5 := \|u(\cdot,t_0)\|_{H^3(\Omega)}
\end{align*}
and we used 
\begin{align*}
|\mathrm{rot}\, L_1 H(x,t_0)|^2 
\le C\left(|\nabla \mathrm{rot}\, H(x,t_0)|^2 
 + |\nabla H(x,t_0)|^2 + |H(x,t_0)|^2\right),                                    \quad x\in \Omega.
\end{align*}
Next, we give the upper bounds of the first term to the fourth term on the right-hand side of 
\eqref{stab3:eq2}, respectively. 

By integration by parts, we estimate
\begin{align*}
\nonumber
&\int_{\Omega} |\partial_t v(x,t_0)|^2 e^{2s\varphi(x,t_0)}dx\\
\nonumber
= &\; \int_{t_0-\delta}^{t_0} \partial_t 
 \left(\int_\Omega |\partial_t v|^2 e^{2s\varphi}dx \right) dt
 + \int_\Omega |\partial_t v(x,t_0-\delta)|^2 e^{2s\varphi(x,t_0-\delta)} dx \\
\nonumber
= &\; \int_{t_0-\delta}^{t_0} \int_\Omega 
\left( 2 (\partial_t v)\cdot (\partial_t^2 v) 
+ 2s(\partial_t \varphi)|\partial_t v|^2 \right)e^{2s\varphi} dxdt\\
\nonumber
& + \int_\Omega |\partial_t v(x,t_0-\delta)|^2 e^{2s\varphi(x,t_0-\delta)} dx \\
\nonumber
\le &\; C\int_Q 
 \left( |\partial_t v| |\partial_t^2 v| + s|\partial_t v|^2 \right)e^{2s\varphi} dxdt 
 + \int_\Omega |\partial_t v(x,t_0-\delta)|^2 e^{2s\varphi(x,t_0-\delta)} dx \\
\nonumber
\le &\; C\int_Q \left(s^{-1}|\partial_t^2 v|^2 
 + s|\partial_t v|^2 \right)e^{2s\varphi} dxdt 
 + Ce^{2s\mu_1}M^2\\
\le &\; C\int_Q \left(s^{-1}|\partial_t^2 v|^2 
 + s^3|\partial_t v|^2 \right)e^{2s\varphi} dxdt 
 + Ce^{2s\mu_1}M^2
\end{align*}
for all $s\ge 1$. 
Here $M>0$ is the constant introduced in the statement of the theorem 
and $\mu_1$ is the same constant introduced in the former subsection, 
see \eqref{con:mu_1}. 

Similarly, we have
\begin{align*}
\int_{\Omega} |H(x,t_0)|^2 e^{2s\varphi(x,t_0)}dx
&\le C\int_Q \left(s^{-1}|\partial_t H|^2 
 + s|H|^2 \right)e^{2s\varphi} dxdt 
 + Ce^{2s\mu_1}M^2 \\
&\le C\int_Q \left(s^4|\partial_t H|^2 
 + s^4|H|^2 \right)e^{2s\varphi} dxdt 
 + Ce^{2s\mu_1}M^2,
\end{align*}
\begin{align*}
\int_{\Omega} |\nabla H(x,t_0)|^2 e^{2s\varphi(x,t_0)}dx
&\le C\int_Q \left(s^{-1}|\nabla (\partial_t H)|^2 
 + s|\nabla H|^2 \right)e^{2s\varphi} dxdt 
 + Ce^{2s\mu_1}M^2 \\
&\le C\int_Q \left(s^2|\nabla (\partial_t H)|^2 
 + s^2|\nabla H|^2 \right)e^{2s\varphi} dxdt 
 + Ce^{2s\mu_1}M^2,
\end{align*}
\begin{align*}
&\int_{\Omega} |\nabla \mathrm{rot}\, H(x,t_0)|^2 e^{2s\varphi(x,t_0)}dx\\
\le &\; C\int_Q \left(s^{-1}|\nabla \mathrm{rot}(\partial_t H)|^2 
 + s|\nabla \mathrm{rot}\, H|^2 \right)e^{2s\varphi} dxdt 
 + Ce^{2s\mu_1}M^2 \\
\le &\; C\int_Q \left(s|\nabla \mathrm{rot}(\partial_t H)|^2 
 + s|\nabla \mathrm{rot}\, H|^2 \right)e^{2s\varphi} dxdt 
 + Ce^{2s\mu_1}M^2
\end{align*}
for all $s\ge 1$. 
Thus, we reach
\begin{align}
\nonumber
&\int_\Omega |\mathrm{rot}\, F(x,t_0)|^2 e^{2s\varphi(x,t_0)}dx \\
\nonumber
\le &\; C\int_Q \left(s^{-1}|\partial_t^2 v|^2 
 + s^3|\partial_t v|^2 \right)e^{2s\varphi} dxdt 
 + C\int_Q \left(s|\nabla \mathrm{rot}(\partial_t H)|^2 
 + s|\nabla \mathrm{rot}\, H|^2 \right)e^{2s\varphi} dxdt \\
\nonumber
&+ C\int_Q \left(s^4|\partial_t H|^2 
 + s^4|H|^2 \right)e^{2s\varphi} dxdt
 + C\int_Q \left(s^2|\nabla (\partial_t H)|^2 
 + s^2|\nabla H|^2 \right)e^{2s\varphi} dxdt \\
\label{stab3:eq4}
& + Ce^{2s\mu_1}M^2 + Ce^{Cs}E_5^2. 
\end{align}

\noindent{\bf Second step: Application of the Carleman estimate} 

For simplicity of the expressions, we introduce the following notations: 
$$
u^{(1)} = \partial_t u, \ 
H^{(1)} = \partial_t H, \ 
p^{(1)} = \partial_t p.
$$
Recall that we take the time derivative of the governing system \eqref{sy:MHD} and obtain
\begin{equation}
\label{sy:MHD1b}
\left\{
\begin{aligned}
&\partial_t u^{(1)} \!-\! \nu\Delta u^{(1)} + (A^{(1)}\!\cdot\!\nabla) u^{(1)} 
 + (u^{(1)}\!\cdot\!\nabla) A^{(2)} + L_1 H^{(1)} + \nabla p^{(1)} = \partial_t F + J_1, \\
&\partial_t H^{(1)} - \kappa\Delta H^{(1)} + (A^{(3)}\cdot\nabla) H^{(1)} 
 + (H^{(1)}\cdot\nabla) A^{(4)} + L_2 u^{(1)} = J_2, \\
& \mathrm{div}\, u^{(1)} = 0,
\end{aligned}
\right.
\end{equation}
where $J_1, J_2$ satisfy
\begin{align*}
|\mathrm{rot}\, J_k| 
\le &\; C\left(|u| + |\nabla u| + |\nabla \mathrm{rot}\, u| 
 + |H| + |\nabla H| + |\nabla \mathrm{rot}\, H|\right),                          \quad k=1,2. 
\end{align*}
We employ Theorem \ref{thm:CE2} to \eqref{sy:MHD} 
and by the same argument in deriving \eqref{stab2:eq6}, we obtain
\begin{align}
\nonumber
&\int_Q \left(s|\nabla \mathrm{rot}\, u|^2 + s|\nabla \mathrm{rot}\, H|^2 
 + s^2|\nabla u|^2 + s^2|\nabla H|^2 + s^4|u|^2 + s^4|H|^2 \right)e^{2s\varphi} dxdt\\
\label{stab3:eq6}
&\le C\int_Q |\mathrm{rot}\, F|^2 e^{2s\varphi} dxdt 
 + Cs^4 e^{Cs}E_6^2 + Cs^4e^{2s\mu_1}M^2
\end{align}
for all sufficiently large $s\ge 1$. Here 
\begin{align*}
E_6 := &\; 
\|\nabla \mathrm{rot}\, u\|_{H^1(I;L^2(\Gamma))} 
+ \|\nabla \mathrm{rot}\, H\|_{H^1(I;L^2(\Gamma))} 
+ \|\nabla u\|_{H^2(I;L^2(\Gamma))} \\
& + \|\nabla H\|_{H^2(I;L^2(\Gamma))} 
+ \|u\|_{H^1(I;L^2(\Gamma))} 
+ \|H\|_{H^1(I;L^2(\Gamma))}.
\end{align*}
Next, we employ Theorem \ref{thm:CE2} to \eqref{sy:MHD1b} and we obtain
\begin{align}
\nonumber
&\int_Q \bigg( \frac{1}{s}|\partial_t \mathrm{rot}\, u^{(1)}|^2 
 + s|\nabla \mathrm{rot}\, H^{(1)}|^2 + s^2|\nabla H^{(1)}|^2 
 + s^3|\mathrm{rot}\, u^{(1)}|^2 + s^4|H^{(1)}|^2
\bigg)e^{2s\varphi} dxdt\\
\nonumber
&\le C\int_Q \left(|\mathrm{rot}(\partial_t F) + \mathrm{rot}\, J_1|^2 
 + |\mathrm{rot}\, J_2|^2\right) e^{2s\varphi} dxdt 
 + Cs^4 e^{Cs}E_6^2 + Cs^4e^{2s\mu_1}M^2\\
\nonumber
&\le C\int_Q |\partial_t \mathrm{rot}\, F|^2 e^{2s\varphi} dxdt 
 + Cs^4 e^{Cs}E_6^2 + Cs^4e^{2s\mu_1}M^2\\
\label{stab3:eq7}
&\quad + C \int_Q \left(|\nabla \mathrm{rot}\, u|^2 + |\nabla \mathrm{rot}\, H|^2 
 + |\nabla u|^2 + |\nabla H|^2 + |u|^2 + |H|^2 \right)e^{2s\varphi} dxdt
\end{align}
for all sufficiently large $s\ge 1$. 
Inserting \eqref{stab3:eq6} into the right-hand side of \eqref{stab3:eq7} yields
\begin{align*}
&\int_Q \bigg( \frac{1}{s}|\partial_t \mathrm{rot}\, u^{(1)}|^2 
 + s|\nabla \mathrm{rot}\, H^{(1)}|^2 + s^2|\nabla H^{(1)}|^2 
 + s^3|\mathrm{rot}\, u^{(1)}|^2 + s^4|H^{(1)}|^2
\bigg)e^{2s\varphi} dxdt\\
&\le C\int_Q \left(|\partial_t \mathrm{rot}\, F|^2 
 + |\mathrm{rot}\, F|^2\right)e^{2s\varphi} dxdt 
 + Cs^4 e^{Cs}E_6^2 + Cs^4e^{2s\mu_1}M^2
\end{align*}
for all sufficiently large $s\ge 1$. 
We combine this inequality with \eqref{stab3:eq6} to obtain
\begin{align*}
&\int_Q \bigg( \frac{1}{s}|\partial_t \mathrm{rot}\, u^{(1)}|^2 
 + s|\nabla \mathrm{rot}\, H^{(1)}|^2 + s|\nabla \mathrm{rot}\, H|^2 
 + s^2|\nabla H^{(1)}|^2 + s^2|\nabla H|^2 \\
&\qquad + s^3|\mathrm{rot}\, u^{(1)}|^2 + s^4|H^{(1)}|^2 + s^4|H|^2
\bigg)e^{2s\varphi} dxdt\\
&\le C\int_Q \left(|\partial_t \mathrm{rot}\, F|^2 
 + |\mathrm{rot}\, F|^2\right)e^{2s\varphi} dxdt 
 + Cs^4 e^{Cs}E_6^2 + Cs^4e^{2s\mu_1}M^2.
\end{align*}
According to the assumption \eqref{con:F2-1}, 
we can estimate the first term on the right-hand side above 
and we further have
\begin{align}
\nonumber
&\int_Q \bigg( \frac{1}{s}|\partial_t \mathrm{rot}\, u^{(1)}|^2 
 + s|\nabla \mathrm{rot}\, H^{(1)}|^2 + s|\nabla \mathrm{rot}\, H|^2 
 + s^2|\nabla H^{(1)}|^2 + s^2|\nabla H|^2 \\
\nonumber
&\qquad + s^3|\mathrm{rot}\, u^{(1)}|^2 + s^4|H^{(1)}|^2 + s^4|H|^2
\bigg)e^{2s\varphi} dxdt\\
\label{stab3:eq8}
&\le C\int_Q |\mathrm{rot}\, F(x,t_0)|^2e^{2s\varphi} dxdt 
 + Cs^4 e^{Cs}E_6^2 + Cs^4e^{2s\mu_1}M^2 
\end{align}
for all sufficiently large $s\ge 1$. 
By noting $v = \mathrm{rot}\, u$, 
we insert \eqref{stab3:eq8} into the right-hand side of \eqref{stab3:eq4} and we obtain
\begin{align}
\nonumber
\int_\Omega |\mathrm{rot}\, F(x,t_0)|^2 e^{2s\varphi(x,t_0)}dx 
\le &\; C\int_Q |\mathrm{rot}\, F(x,t_0)|^2 e^{2s\varphi} dxdt \\
\label{stab3:eq9}
& + Cs^4e^{2s\mu_1}M^2 + Cs^4e^{Cs}D_3^2
\end{align}
for all sufficiently large $s\ge 1$. 
Here we noted that $E_5^2 + E_6^2 \le C D_3^2$ where 
$D_3$ is the norm of the data that introduced in the theorem. 

\noindent {\bf Third Step: Completion of the proof}

As we mentioned in the fourth step in Subsection \ref{subsec:proof1}, 
thanks to the choice of the weight function $\varphi$, we find
\begin{align*}
\int_Q |\mathrm{rot}\, F(x,t_0)|^2 e^{2s\varphi} dxdt 
= o(1)\int_\Omega |\mathrm{rot}\, F(x,t_0)|^2 e^{2s\varphi(x,t_0)}dx
\end{align*}
as $s\to \infty$. 
Thus, by taking $s$ large enough, we can absorb the first term on the right-hand side of 
\eqref{stab3:eq9} into the left-hand side. This yields
\begin{align}
\label{stab3:eq10}
\int_\Omega |\mathrm{rot}\, F(x,t_0)|^2 e^{2s\varphi(x,t_0)} dx 
\le &\; Cs^4e^{2s\mu_1}M^2 + Cs^4e^{Cs}D_3^2
\end{align}
for all $s$ large enough. 
With the same argument and the same constant $\mu_0>0$ introduced in the last step in 
Subsection \ref{subsec:proof2}, 
we have the lower bound of the left-hand side of \eqref{stab3:eq10}:
\begin{align*}
\int_\Omega |\mathrm{rot}\, F(x,t_0)|^2 e^{2s\varphi(x,t_0)} dx 
\ge e^{2s\mu_0} \|\mathrm{rot}\, F(\cdot,t_0)\|_{L^2(\Omega_0)}^2,  
\end{align*} 
which implies
\begin{align}
\label{stab3:eq11}
\|\mathrm{rot}\, F(\cdot,t_0)\|_{L^2(\Omega_0)}^2
\le Cs^4e^{-2s(\mu_0-\mu_1)}M^2 + Cs^4e^{Cs}D_3^2
\end{align}
for all $s$ large enough. 
Since we have $\mu_0 > \mu_1$ by taking $\beta>0$ sufficiently large, 
we end up with the stability of H\"older type by minimizing 
the right-hand side of \eqref{stab3:eq11} with respect to $s$. 
%
%
\section{Concluding remarks}
\subsection{Linearization of the MHD equations}
Let $(u_i,p_i,H_i)$ satisfy the following MHD equations
$$
\left\{
\begin{aligned}
& \ \partial_t u_i - \nu\Delta u_i + (u_i\cdot\nabla) u_i 
 - \mu\, \mathrm{rot}\, H_i\times H_i + \nabla p_i = F_i, \\
& \ \partial_t H_i - \frac{1}{\sigma\mu}\Delta H_i - \mathrm{rot}(u_i\times H_i) = 0, \\
& \ \mathrm{div}\, u_i = 0, \quad \mathrm{div}\, H_i = 0
\end{aligned}   
\right.                                                 
$$
for $i=1,2$. Then by taking the difference and setting 
$u=u_1-u_2$, $p=p_1-p_2$, $H=H_1-H_2$ and $F=F_1-F_2$, 
we have
$$
\left\{
\begin{aligned}
& \ \partial_t u - \nu\Delta u + (u_1\cdot\nabla) u 
 + (u\cdot\nabla) u_2 - \mu (H_1\cdot\nabla) H - \mu (H\cdot\nabla) H_2 \\
&\quad + \frac{\mu}{2} \nabla ((H_1+H_2)\cdot H) + \nabla p = F, \\
& \ \partial_t H - \frac{1}{\sigma\mu}\Delta H + (u_1\cdot\nabla) H 
 - (H\cdot\nabla) u_2 - (H_1\cdot\nabla) u + (u\cdot\nabla) H_2 = 0, \\
& \ \mathrm{div}\, u = 0, \quad \mathrm{div}\, H = 0,
\end{aligned}   
\right.                                                 
$$
which corresponds to the linearized MHD equations \eqref{sy:MHD} with \eqref{def:L1L2}. 
For example, we let $A^{(1)} = u_1$, $A^{(2)} = u_2$, \ldots. 

Therefore, in order to apply the main results in this article to the (nonlinear) MHD equations, 
we need some a priori regularity of the solutions $(u_i,p_i,H_i)$, $i=1,2$. 
Although the existence of the solutions to the initial boundary value problems 
for the MHD equations is a serious issue with the regularity, 
here we focus on the inverse problems and we assume the existence of such solutions 
$(u_i,p_i,H_i)$, $i=1,2$ with sufficient regularity.

\subsection{Global Lipschitz stability}
As for the Carleman estimates for the parabolic equations including the Navier-Stokes equations 
and our MHD equations, one can have Carleman estimates according to each of the weight functions 
in the following forms:
\begin{equation}
\label{weight:reg}
\mbox{(i)} \qquad \varphi(x,t) := e^{\lambda (d(x)-\beta(t-t_0)^2)},
\end{equation}
which is regular in $x,t$.
\begin{equation}
\label{weight:sin}
\mbox{(ii)} \qquad \alpha(x,t) := \frac{e^{\lambda d(x)} 
- e^{2\lambda \|d\|_{C(\overline{\Omega})}}}{\ell(t)},
\end{equation}
which is singular at $t=0,T$.

Usually one uses the Carleman estimate with the regular weight function \eqref{weight:reg} 
to prove a local H\"older stability estimate for an inverse problem. 
On the other hand, one can apply the Carleman estimate with the singular weight function 
\eqref{weight:sin} to prove a global Lipschitz stability for an inverse problem provided that 
the boundary condition on $\partial\Omega \times (0,T)$ is given. 

In this article, it is clear that we employed the Carleman estimates with \eqref{weight:reg} 
and proved the H\"older stability estimates for some inverse source problems. 
If we assume suitable boundary conditions, for example, 
$$
u(x,t) = 0,\quad 
H(x,t) = 0,\quad          (x,t)\in \partial\Omega\times (0,T),
$$
in additional to \eqref{sy:MHD}, we can easily combine the arguments in 
Choulli, Imanuvilov, Puel and Yamamoto \cite{CIPY13} or Imanuvilov and Yamamoto \cite{IY21}, 
in which the authors discussed similar inverse source problems for the linearized Naiver-Stokes 
equations by the Carleman estimates with \eqref{weight:sin}, 
with what we have done in this article to derive the global Lipschitz stability estimates for 
the inverse source problems. 

However, the boundary condition on the magnetic field intensity 
$H(x,t) = 0$, $(x,t)\in \partial\Omega\times (0,T)$ is not physical. 
As a possible solution, instead of the boundary conditions, we measure $u$ and $H$ 
on the full boundary. 
Then immediately we can prove the related global Lipschitz stability which is parallel to 
Theorem \ref{thm:stab1} as follows:
\begin{align*}
&\|F(\cdot,t_0)\|_{L^2(\Omega)} 
\le C( \|u\|_{H^2(I;H^1(\partial\Omega))} 
 + \|\partial_n u\|_{H^1(I;L^2(\partial\Omega))} 
 + \|H\|_{H^2(I;H^1(\partial\Omega))} \\
&\qquad + \|\partial_n H\|_{H^1(I;L^2(\partial\Omega))} 
 + \|p\|_{H^1(I;H^{\frac12}(\partial\Omega))} 
 + \|u(\cdot,t_0)\|_{H^2(\Omega)} 
 + \|\nabla p(\cdot,t_0)\|_{L^2(\Omega)}).
\end{align*}
Comparing to the inverse source problem for the parabolic equations, 
the data on the right-hand side above are too many since we need 
additionally $\partial_n u,\partial_n H$ on the boundary. 
The global Lipschitz stability for the inverse source problems for the MHD equations by less data 
still remains open. 
\section*{Acknowledgments}
The first author thanks the Leading Graduate Course 
for Frontiers of Mathematical Sciences and Physics 
(FMSP, The University of Tokyo), 
and was supported by Grant-in-Aid for JSPS Fellows 20F20319 of 
Japan Society for the Promotion of Science (JSPS). 
The second author was supported by Grant-in-Aid for Scientific Research (A) 20H00117 of JSPS, 
the National Natural Science Foundation of China (Nos. 11771270, 91730303) 
and the RUDN University Strategic Academic Leadership Program. 


\vspace{0.2cm}

E-mail address: huangxc@ms.u-tokyo.ac.jp

E-mail address: myama@ms.u-tokyo.ac.jp
\end{document}